\documentclass[%
 aip,
 amsmath,amssymb,
 reprint,%
]{revtex4-1}
\usepackage{graphicx}% Include figure files
\usepackage{dcolumn}% Align table columns on decimal point
\usepackage{bm}% bold math
\usepackage[cal=cm]{mathalfa}
\usepackage[utf8]{inputenc}
\usepackage[T1]{fontenc}
\usepackage{mathptmx}
\usepackage{etoolbox}
\usepackage{algorithm} 
\usepackage{subfigure}
\usepackage{algpseudocode} 
\usepackage{amsthm}
\usepackage{multirow}
\newtheorem{remark}{Remark}
\newtheorem{lemma}{Lemma}[section]
\newtheorem{thm}{\bf Theorem}[section]

\begin{document}

\preprint{AIP/123-QED}

\title[]{Runge-Kutta Random Feature Method for Solving  Multiphase Flow Problems of Cells}
% Force line breaks with \\
\author{Yangtao Deng}

 % \altaffiliation{School of Mathematics, Sichuan University, Chengdu, China}%Lines break automatically or can be forced with \\
% \author{Qiaolin He}%
%  \email{qlhejenny@scu.edu.cn}
% \affiliation{ 
% School of Mathematics, Sichuan University, Chengdu, China
% }%

\author{Qiaolin He}
\email[Corresponding author: ]{qlhejenny@scu.edu.cn}
\affiliation{\textsuperscript{1} 
School of Mathematics, Sichuan University, Chengdu, China
}%

\date{\today}% It is always \today, today,

\begin{abstract}
Cell collective migration plays a crucial role in a variety of physiological processes. In this work, we propose the Runge-Kutta random feature method to solve the nonlinear and strongly coupled multiphase flow problems of cells, in which  the random feature method in space and the explicit Runge-Kutta  method in time are utilized. Experiments indicate that this algorithm can effectively deal with time-dependent partial differential equations  with strong nonlinearity, and achieve high accuracy both in space and time.  Moreover, in order to improve computational efficiency and save computational resources, we choose to implement parallelization and non-automatic differentiation strategies in our simulations. We also provide error estimates for the Runge-Kutta random feature method, and a series of numerical experiments are shown to validate our method.
\end{abstract}

\maketitle

% \begin{quotation}
% The ``lead paragraph'' is encapsulated with the \LaTeX\ 
% \verb+quotation+ environment and is formatted as a single paragraph before the first section heading. 
% (The \verb+quotation+ environment reverts to its usual meaning after the first sectioning command.) 
% Note that numbered references are allowed in the lead paragraph.
% %
% The lead paragraph will only be found in an article being prepared for the journal \textit{Chaos}.
% \end{quotation}

\section{Introduction}\label{sec01}

Cell collective migration plays a vital role in essential physiological processes. Research work about how mechanical forces interact with collective cell behavior at the tissue level \cite{ladoux2017mechanobiology} has garnered growing interest across various fields. A phase field system \eqref{eq:multiphasepde} is used to model cells as active deformable droplets in two dimensions in reference \cite{mueller2019emergence}. Although finite difference methods \cite{leveque2007finite}, finite volume methods \cite{moukalled2016finite}, and finite element methods 
\cite{thomee2007galerkin} can be employed to solve this kind of partial differential equations (PDEs) and great success in theory and application have been achieved, these methods still face challenges such as complex geometric shapes, mesh generation, and lack of generalization. These approaches generally involve generating the grid of points or elements, discretizing differential and integral expressions, and solving the resulting discrete equations. However, the complexity of the computational domain can make mesh generation very difficult  and the problem discretization could create a bias between the mathematical nature of the PDEs and its approximating model. Moreover, obtaining numerical values outside the computational domain requires regenerating the grid and recalculating, which can waste a lot of computational resources and time.

On the other hand, the success of deep learning in computer vision and natural language processing \cite{goodfellow2016deep} has attracted widespread attention in the field of scientific computing. As a special class of functions, neural networks have been proven to be universal approximators of continuous functions \cite{cybenko1989approximation}. Many researchers are committed to using neural networks to solve ordinary differential equations (ODEs) and PDEs \cite{han2017deep,han2018solving,yu2018deep,sirignano2018dgm,zang2020weak,raissi2019physics,weinan2021algorithms}. %Since PDE solutions can be defined in variational form (if they exist), strong form, and weak form, 
Due to the different formulations of PDEs, Deep Ritz methods (DRM) \cite{yu2018deep}, Deep Galerkin methods (DGM) \cite{sirignano2018dgm}, physics-informed neural networks (PINN) \cite{raissi2019physics}, and weak adversarial networks (WAN) \cite{zang2020weak} have been proposed based on loss (objective) functions of variational, strong, and weak forms, respectively. Deep learning-based algorithms have become general methods to solve various PDEs without any type of mesh generation. However, there is a lack of reliable error control in existing deep learning methods. For example, in the absence of the exact solution, as the number of parameters increases, the numerical approximation provided by deep learning methods does not show a clear trend of convergence.

Recent studies have shown that a specific type of neural network - extreme learning machine (ELM) \cite{huang2006extreme} or random features \cite{neal2012bayesian,rahimi2007random} can be used to solve ODEs and PDEs. ELM is a feedforward neural network with a single hidden layer, which randomly assigns weights and biases between the input layer and the hidden layer, and determines the weights of the output layer by analytical methods, see \cite{huang2015trends}. Due to the characteristics of this architecture, the weights and biases of the hidden layer do not need to be learned. This makes ELM faster than typical deep neural networks, because optimization may lead to extremely slow learning speeds in deep neural networks. Overall, randomized neural networks improve the efficiency of learning tasks while maintaining high accuracy in numerical solutions. ELM is a mesh-free method, so it can easily handle the PDEs in complex geometric domains \cite{chen2022bridging,chen2023random}. In addition, ELM is a universal approximator \cite{huang2006extreme,huang2006universal}, so it has the potential to represent any continuous function.
In reference \cite{dwivedi2020physics}, ELM is developed to differential problems, where the authors implemented a model called physics-informed extreme learning machine (PIELM), and subsequently ELM has been successfully applied to solve ODEs and PDEs \cite{chen2022bridging,chen2023random,calabro2021extreme,dong2021local,dwivedi2020physics,fabiani2021numerical,yang2018novel}.

%The first use of ELM in differential problems can be traced back to \cite{dwivedi2020physics}, 

In reference \cite{dong2021local}, a special type of partition of unity (PoU), called domain decomposition, is combined with the ELM to approximate the solution of PDEs, and a block time-marching strategy is proposed for long-time simulations, while a nonlinear least squares method is used to solve nonlinear differential equations. However, as the simulation time increases, the error will grow exponentially quickly \cite{chen2023random}, in addition, when dealing with nonlinear differential equations that exhibit strong nonlinearity (e.g.\eqref{eq:multiphasepde}), the nonlinear least squares method is difficult to converge to satisfactory results. In reference \cite{chen2022bridging}, combining PoU and random feature functions, the random feature method (RFM) is proposed to solve static linear PDEs with complex geometric shapes. In reference \cite{chen2023random}, based on RFM, the  space-time random feature method (ST-RFM) is proposed to solve linear time-dependent PDEs using random feature functions  that depend on both spatial and temporal variables (STC), or the product of two random feature functions that depend on spatial and temporal variables separately (SoV). In addition, an automatic rescaling strategy is proposed to balance the weights of equations and boundary conditions. However, since the time variable is used as the input of the random feature functions, the ST-RFM requires a large amount of computational resources in the process, and can only solve linear differential equations.

In this work, in order to solve the strongly nonlinear coupled multiphase flow problem of cells \eqref{eq:multiphasepde} with less computational resources while maintaining high numerical accuracy, we choose to combine traditional numerical methods with deep learning-based methods. We use the RFM in space and the explicit Runge-Kutta (RK) method \cite{butcher1996history} in time, which is called the Runge-Kutta random feature method (RK-RFM). This new type of algorithm can achieve high accuracy in space and time. Besides, it is mesh-free, making them easy to use even in complex geometric settings. Moreover, this method can effectively solve PDEs with strong nonlinearity. At the same time,  we choose to use the same set of random feature functions for different cells and manually derive the necessary neural network derivative functions before starting the calculations, instead of using automatic differentiation, which greatly accelerates computation and saves computational resources. We also provide the error estimates for the RK-RFM, and show some numerical experiments. 

 The rest of this article is organized as follows. In Section \ref{sec02}, we briefly introduce RFM and the multiphase flow model of cells. In Section \ref{sec03}, we present the RK-RFM and give  the error estimates of method. In Section \ref{sec04}, we show some numerical experiments to verify our theoretical results and the application of the RK-RFM in the multiphase flow problem of cells. The conclusions and remarks are given in Section \ref{sec05}.

\section{Random Feature Method and Multiphase Flow Problem of Cells}
\label{sec02}

\subsection{Random Feature Method for Static Problems}
\label{sec02:sec01}
% Let $\boldsymbol{x} \in \Omega \subset \mathbb{R}^{d_x}$, where $d_x \in \mathbb{N}^{+}$is the dimension of $\boldsymbol{x}$, and let $d_u \in \mathbb{N}^{+}$be the dimension of the output. 
Consider the following  linear boundary-value problem
\begin{equation}
	\begin{aligned}
\begin{cases}\mathcal{L} \boldsymbol{\phi}(\boldsymbol{x})=\boldsymbol{f}(\boldsymbol{x}), & \boldsymbol{x} \in \Omega, \\ \mathcal{B} \boldsymbol{\phi}(\boldsymbol{x})=\boldsymbol{g}(\boldsymbol{x}), & \boldsymbol{x} \in \partial \Omega,\end{cases}
	\end{aligned}\label{eq:spde}
\end{equation}
where $\mathcal{L}$ and $\mathcal{B}$ are linear differential and linear boundary operators, respectively. We use $d$ and $d_{\phi}$ to denote the dimensions of $\boldsymbol{x}$ and $\boldsymbol{\phi} = (\phi_1, \phi_2,...,\phi_{d_\phi})^T$, respectively.

%The RFM \cite{chen2022bridging} primarily includes the following steps. 
In RFM \cite{chen2022bridging}, the domain $\Omega$ is divided into $M_p$ non-overlapping subdomains ${\Omega_n}$ with $\boldsymbol{x}_n$ as the centers.  We have $\Omega =\cup_{n=1}^{M_p} \Omega_n$. For each $\Omega_n$, the linear transformation is employed in RFM
\begin{equation}
	\begin{aligned}
\tilde{\boldsymbol{x}}=\frac{1}{\boldsymbol{r}_n}\left(\boldsymbol{x}-\boldsymbol{x}_n\right), \quad n=1, \cdots, M_p,
\end{aligned}\label{eq:normalize}
\end{equation}
to  map $\Omega_n$ into $[-1,1]^d$, where $\boldsymbol{r}_n\in R^d$ is the radius of $\Omega_n$. A PoU function $\psi_n$ is defined such that $\operatorname{supp}(\psi_n) = \Omega_n$. Two commonly used PoU functions for $d=1$ are
% \begin{equation}
% 	\begin{aligned}
% & \psi(\boldsymbol{x})=\mathbb{I}_{[-1,1]}(\tilde{\boldsymbol{x}}),
% \end{aligned}\label{eq:psi1}
% \end{equation}
% \begin{equation}
% 	\begin{aligned}
% & \psi(\boldsymbol{x})=\mathbb{I}_{\left[-\frac{5}{4},-\frac{3}{4}\right]}(\tilde{\boldsymbol{x}}) \frac{1+\sin (2 \pi \tilde{\boldsymbol{x}})}{2}+\mathbb{I}_{\left[-\frac{3}{4}, \frac{3}{4}\right]}(\tilde{\boldsymbol{x}})\\
% &+\mathbb{I}_{\left[\frac{3}{4}, \frac{5}{4}\right]}(\tilde{\boldsymbol{x}}) \frac{1-\sin (2 \pi \tilde{\boldsymbol{x}})}{2}.
% \end{aligned}\label{eq:psi2}
% \end{equation}
\begin{eqnarray}
%	\begin{aligned}
& &\psi(\boldsymbol{x})  = \mathbb{I}_{[-1,1]}(\tilde{\boldsymbol{x}}) \label{eq:psi1} %\\ 
  %& & \nonumber \\
\end{eqnarray}
and 
\begin{eqnarray}
& &  \psi(\boldsymbol{x})  = \mathbb{I}_{\left[-\frac{5}{4},-\frac{3}{4}\right]}(\tilde{\boldsymbol{x}}) \frac{1+\sin (2 \pi \tilde{\boldsymbol{x}})}{2}+\mathbb{I}_{\left[-\frac{3}{4}, \frac{3}{4}\right]}(\tilde{\boldsymbol{x}})\label{eq:psi2}\\
& & +\mathbb{I}_{\left[\frac{3}{4}, \frac{5}{4}\right]}(\tilde{\boldsymbol{x}}) \frac{1-\sin (2 \pi \tilde{\boldsymbol{x}})}{2}, \nonumber
%\end{aligned}
\end{eqnarray}
where $\mathbb{I}_{[a,b]}(\boldsymbol{x})=1, \boldsymbol{x}\in[a,b]$ and $a\le b$. When $d>1$, the PoU function $\psi(\boldsymbol{x})$ is defined as $\psi(\boldsymbol{x})=\Pi_{i=1}^{d} \psi\left(x_i\right)$.

Next, random feature functions $\varphi_{n j}$ on $\Omega_n$ is formulated by a two-layer neural network
\begin{equation}
	\begin{aligned}
\varphi_{n j}(\boldsymbol{x})=\sigma\left(\boldsymbol{W}_{n j}\cdot \tilde{\boldsymbol{x}}+b_{n j}\right), \quad j=1,2, \cdots, J_n,
\end{aligned}\label{eq:basicfunction}
\end{equation}
where $\sigma$ is the nonlinear activation function, the $\boldsymbol{W}_{n j}$ and $b_{n j}$ are randomly chosen through uniform distribution $\mathbb{U}\left(-R_m, R_m\right)$ but fixed, the $R_m$ controls the magnitude of the parameters and $J_n$ is the number of random feature functions. In RFM, the approximate solution is formed by the linear combination of these random feature functions along with the PoU functions
\begin{equation}
	\begin{aligned}
\tilde{\boldsymbol{\phi}}(\boldsymbol{x})=&\left(\sum_{n=1}^{M_p} \psi_n(\boldsymbol{x}) \sum_{j=1}^{J_n} u_{n j}^{1} \varphi_{n j}^{1}(\boldsymbol{x}),\cdots,\right.\\
&\left.\sum_{n=1}^{M_p} \psi_n(\boldsymbol{x}) \sum_{j=1}^{J_n} u_{n j}^{d_{\phi}} \varphi_{n j}^{d_{\phi}}(\boldsymbol{x})\right)^T,
\end{aligned}\label{eq:approximatesolution}
\end{equation}
where $u_{n j}^{i} \in \mathbb{R}$ is unknowns to be determined, and $M=M_p J_n$ denotes the degree of freedom.

Then, RFM minimizes the loss function
\begin{equation}
	\begin{aligned}
\operatorname{Loss}\left(\left\{u_{n j}^{i}\right\}\right)=&\sum_{n=1}^{M_p}\left(\sum_{q=1}^Q\left\|\boldsymbol{\lambda}_{n, q}\left(\mathcal{L} \tilde{\boldsymbol{\phi}}\left(\boldsymbol{x}_{q}^{n}\right)-\boldsymbol{f}\left(\boldsymbol{x}_{q}^{n}\right)\right)\right\|_2^2\right)\\
&+\sum_{n=1}^{M_p}\left(\sum_{\boldsymbol{x}_{q}^{n} \in \partial \Omega}\left\|\boldsymbol{\lambda}_{n, b}\left(\mathcal{B}\tilde{\boldsymbol{\phi}}\left(\boldsymbol{x}_{q}^{n}\right)-\boldsymbol{g}\left(\boldsymbol{x}_{q}^{n}\right)\right)\right\|_2^2\right)
\end{aligned}\label{eq:loss}
\end{equation}
 though the linear least-squares method.
When employing the PoU function $\psi$ defined in \eqref{eq:psi1}, it is necessary to incorporate regularization terms into the loss function \eqref{eq:loss} to ensure continuity between neighboring subdomains. Conversely, when utilizing $\psi$ defined in \eqref{eq:psi2}, regularization terms are not necessary. 
To find the optimal set of parameters $\left\{u_{n j}^{i}\right\}$, RFM samples $Q$ collocation points $\left\{\boldsymbol{x}_{q}^{n}\right\}$ in each $\Omega_n$ and computes the rescaling parameters  $\boldsymbol{\lambda}_{n, q}=diag(\lambda_{n, q}^{1},\cdots,\lambda_{n, q}^{d_{\phi}})$ and $\boldsymbol{\lambda}_{n, b}=diag(\lambda_{n, b}^{1},\cdots,\lambda_{n, b}^{d_{\phi}})$. Specifically, the rescaling parameters are calculated using the following method
\begin{equation}
\begin{aligned} 
& \lambda_{n, q}^{i}=\frac{c}{\mathop{\max}\limits_{1\leq j \leq J_n}\left|\mathcal{L}\left(\psi_n(\boldsymbol{x}_{q}^{n})\varphi_{n j}^{i}(\boldsymbol{x}_{q}^{n})  \right)\right|}, \\
&\quad q=1, \cdots, Q, \ n=1, \cdots, M_p,  \ i=1, \cdots, d_{\phi}, \\
&\lambda_{n, b}^{i}=\frac{c}{\mathop{\max}\limits_{1\leq j \leq J_n}\left|\mathcal{B}\left(\psi_n(\boldsymbol{x}_{q}^{n})\varphi_{n j}^{i}(\boldsymbol{x}_{q}^{n})  \right)\right|}, \\
&\quad \boldsymbol{x}_{q}^{n}\in{\partial\Omega}, \ n=1, \cdots, M_p, \ i=1, \cdots, d_{\phi},
\end{aligned}\label{eq:RFMlambda}
\end{equation}
where $c>0$ is a constant, see more details in Appendix \ref{app:A}. The numerical result is formulated by using equation \eqref{eq:approximatesolution}.

\subsection{Multiphase Flow Problem of Cells}
\label{sec02:sec03}
We consider the following time-dependent multiphase flow problem, which is a phase-field problem used to model cells as active deformable droplets as in reference \cite{mueller2019emergence}
\begin{equation}
        \left\{
	\begin{aligned}
 &\partial_t \phi_i=-v_i \nabla \phi_i-\frac{\delta \mathcal{F}}{\delta \phi_i}=F(\boldsymbol{x}, t, \boldsymbol{\phi}(\boldsymbol{x}, t)), \\
 &\boldsymbol{x}, t \in  \Omega \times(0, T],\quad  i=1, \ldots, d_{\phi},\\
% \left\{\begin{array}{lr}
% \frac{\partial\boldsymbol{\phi}}{\partial t}=\mathcal{L} \boldsymbol{\phi}+\boldsymbol{f}(\boldsymbol{x}, t)=\mathcal{F}(\boldsymbol{x}, t, \boldsymbol{\phi}(\boldsymbol{x}, t)), & \boldsymbol{x}, t \in \Omega \times[0, T], \\
&\mathcal{B} \boldsymbol{\phi}(\boldsymbol{x}, t)=\boldsymbol{g}(\boldsymbol{x}, t), \quad \boldsymbol{x}, t \in \partial \Omega \times[0, T], \\
&\boldsymbol{\phi}(\boldsymbol{x}, 0)=\boldsymbol{h}(\boldsymbol{x}), \quad \boldsymbol{x} \in \Omega,
% \end{array}\right.
\end{aligned}\label{eq:multiphasepde}
\right.
\end{equation}
% where $\mathcal{L}$ is a differential operator acting only on the spatial variables and $\mathcal{B}$ is a linear boundary operator. Here, $\mathcal{L}$ does not need to be linear.
where $\mathcal{B}$ is a linear boundary operator and $\phi_i$ represents the i-th cell, with $\phi_i=1$ and $\phi_i=0$ indicating the interior and exterior of the cell, respectively. The cell boundary is defined to locate at $\phi_i=\frac{1}{2}$, and  $\mathcal{F}$ is a free energy  with $\mathcal{F}=\mathcal{F}_{\mathrm{CH}}+\mathcal{F}_{\text {area }}+\mathcal{F}_{\text {rep }}$,
where
\begin{equation}
\begin{aligned}
\mathcal{F}_{\mathrm{CH}} & =\sum_i \frac{\gamma}{\lambda} \int \left\{4 \phi_i^2\left(1-\phi_i\right)^2+\lambda^2\left(\nabla \phi_i\right)^2\right\} d \boldsymbol{x}, \\
\mathcal{F}_{\text {area }} & =\sum_i \mu\left(1-\frac{1}{\pi R^2} \int  \phi_i^2 d \boldsymbol{x}\right)^2, \\
\mathcal{F}_{\text {rep }} & =\sum_i \sum_{j \neq i} \frac{\kappa}{\lambda} \int \phi_i^2 \phi_j^2 d \boldsymbol{x}.
\end{aligned}\label{eq:F}
\end{equation}

The Cahn-Hilliard free energy, denoted as $\mathcal{F}_{\mathrm{CH}}$, is instrumental in maintaining the stability of the cell interfaces. The term $\mathcal{F}_{\text {area }}$ acts as a soft regulation for the area occupied by each cell, aiming for an area close to 
$\pi R^2$, with $R$ representing the radius of the cell. This ensures that the cells have the capacity to be compressed.
Additionally, the repulsion term, $\mathcal{F}_{\text {rep }}$, imposes a penalty on areas where two cells are in overlap. The normalization process has been carefully selected to ensure that the interfaces have a width of $\lambda$ at equilibrium, and that the essential characteristics of the cells remain largely unaffected even when $\lambda$ is adjusted. The parameters 
$\gamma, \mu,$ and $\kappa$ are pivotal in determining the time it takes for shape deformations, area adjustments, and repulsive forces to reach equilibrium, respectively.

Then, taking variation and we have 
\begin{equation}
\begin{aligned}
\frac{\delta \mathcal{F}_{\mathrm{CH}}}{\delta \phi_i}&=\frac{8 \gamma}{\lambda} \phi_i\left(1-\phi_i\right)\left(1-2 \phi_i\right)-2 \gamma \lambda \Delta \phi_i,\\
\frac{\delta \mathcal{F}_{\text {area }}}{\delta \phi_i}&=-\frac{4 \mu}{\pi R^2} \phi_i\left(1-\frac{1}{\pi R^2} \int  \phi_i^2 d \boldsymbol{x}\right),\\
\frac{\delta \mathcal{F}_{\text {rep }}}{\delta \phi_i}&=\frac{2 \kappa}{\lambda} \sum_{k \neq i} \phi_k^2 \phi_i.
\end{aligned}\label{eq:partialFpartialphi}
\end{equation}
The $v_i$ is the total velocity of the i-th cell, given by the following force balance equation 
 \begin{equation}
\begin{aligned}
\xi v_i=F_i^{\text {int }},
\end{aligned}\label{eq:forcebalance}
\end{equation}
where $\xi$ is a substrate friction coefficient and $F_i^{\text {int }}$ is the total force acting on the interface of i-th cell. These microscopic interface forces in terms of a macroscopic tissue stress tensor $\sigma_{\text {tissue }}$  are defined as

 \begin{equation}
\begin{aligned}
F_i^{\text {int }}=\int  \phi_i \nabla\cdot \sigma_{\text {tissue }} d \boldsymbol{x}=-\int  \sigma_{\text {tissue }} \cdot \nabla \phi_i d \boldsymbol{x},
\end{aligned}\label{eq:Fint}
\end{equation}
%The first expression is the integral of the local force $\nabla\cdot\sigma_{\text {tissue }}$ weighted by the phase field $\phi_i$, while the second is the integral of the force exerted by the stress tensor on the vector $-\nabla \phi_i$ normal to the interface and pointing outwards. 
where the tissue stress tensor $\sigma_{\text {tissue }}$  is decomposed into passive and active stresses
\begin{equation}
\begin{aligned}
\sigma_{\text {tissue }}=-P \mathbb{I}-\zeta Q,
\end{aligned}\label{eq:sigmatissue}
\end{equation}
where $\zeta$ represents activity. The fields $P$ and $Q$ are the tissue pressure and tissue nematic tensor, defined as follows 
\begin{equation}
\begin{aligned}
P=\sum_i\left(\frac{\delta \mathcal{F}_{\text {rep }}}{\delta \phi_i}-\frac{\delta \mathcal{F}_{\mathrm{CH}}}{\delta \phi_i}-\frac{\delta \mathcal{F}_{\text {area }}}{\delta \phi_i}\right),
\end{aligned}\label{eq:P}
\end{equation}
and
\begin{equation}
\begin{aligned}
Q=\sum_i \phi_i S_i,
\end{aligned}\label{eq:Q}
\end{equation}
where $S_i$ is the deformation tensor of i-th cell defined as the traceless part of $-\int \left(\nabla \phi_i\right)^{\top} \nabla \phi_i d \boldsymbol{x}$. 
Specifically, $S_i \equiv S\left(\phi_i\right)$, with
\begin{equation}
\begin{aligned}
&S(\boldsymbol{\phi})\\
&=\left(\begin{array}{cc}
S_{11} & S_{12} \\
S_{12} & -S_{11}
\end{array}\right)\\
&=\int \left(\begin{array}{cc}
\frac{1}{2}\left(\left(\partial_y \boldsymbol{\phi}\right)^2-\left(\partial_x \boldsymbol{\phi}\right)^2\right) & -\left(\partial_x \boldsymbol{\phi}\right)\left(\partial_y \boldsymbol{\phi}\right) \\
-\left(\partial_x \boldsymbol{\phi}\right)\left(\partial_y \boldsymbol{\phi}\right) & \frac{1}{2}\left(\left(\partial_x \boldsymbol{\phi}\right)^2-\left(\partial_y \boldsymbol{\phi}\right)^2\right) 
\end{array}\right)d \boldsymbol{x}.
\end{aligned}\label{eq:Sphi}
\end{equation}

% Due to that the time variable $t$ is used as the input to the random feature functions \eqref{eq:STC} and \eqref{eq:SoV}, when $d=2$, the ST-RFM consumes a large amount of computational resources during the calculation process and can only solve linear differential equations. \cite{dong2021local} proposes a block time-marching strategy and uses nonlinear least squares method to solve nonlinear differential equations. The block time-marching strategy is to solve \eqref{eq:tdpde} in each time subdomain sequentially and use the prediction at the terminal time of the previous subdomain as the initial condition for the current subdomain. However, The block time-marching strategy still consumes a lot of computational resources in the calculation process, and it has been proven in \cite{chen2023random} that the error introduced by this strategy grows exponentially with time. Furthermore, when dealing with nonlinear differential equations that exhibit strong nonlinearity (e.g.\eqref{eq:multiphasepde}), the nonlinear least squares method struggles to converge to a satisfactory result. In order to solve the coupled multiphase flow problems of cells \eqref{eq:multiphasepde} with strong nonlinearity under limited computational resources while maintaining high numerical accuracy, we choose to use the RFM in space and an explicit Runge-Kutta method in time, which is named the Runge-Kutta Random Feature Method.

\section{Runge-Kutta Random Feature Method and Error Estimation}
\label{sec03}
\subsection{Runge-Kutta Random Feature Method}
\label{sec03:sec01}
 We now begin to introduce the RK-RFM. Given a partition of $[0,T]: 0 = t_{0}< t_{1}<...<t_{K}=T$ and for $k = 0,...,K-1$, the explicit scheme \cite{butcher1996history} % $p$th-order RK method \cite{butcher1996history} 
 is written as
\begin{equation}\label{eq:porderrk}
	\left\{ 
	\begin{aligned}
        \boldsymbol{D}_1(\boldsymbol{x}) &= F(\boldsymbol{x}, t_{k}, \boldsymbol{\phi}_{t_{k}}(\boldsymbol{x})),\\
        \boldsymbol{D}_2(\boldsymbol{x}) &= F(\boldsymbol{x}, t_{k}+c_2\Delta t, \boldsymbol{\phi}_{t_{k}}(\boldsymbol{x})+a_{21}\Delta t\boldsymbol{D}_1(\boldsymbol{x})),\\
        & \quad \vdots\\
        \boldsymbol{D}_p(\boldsymbol{x}) &= F(\boldsymbol{x}, t_{k}+c_p\Delta t, \boldsymbol{\phi}_{t_{k}}(\boldsymbol{x})+\Delta t\sum_{i=1}^{p-1}a_{p i}\boldsymbol{D}_i(\boldsymbol{x})),\\
        \boldsymbol{\phi}_{t_{k+1}}(\boldsymbol{x}) &= \boldsymbol{\phi}_{t_{k}}(\boldsymbol{x})+\Delta t\sum_{i=1}^{p}b_{i}\boldsymbol{D}_i(\boldsymbol{x})\\
        & \equiv \boldsymbol{\phi}_{t_{k}}(\boldsymbol{x})+\Delta t H(\boldsymbol{x}, t_{k}, \boldsymbol{\phi}_{t_{k}}(\boldsymbol{x}),\Delta t),
	\end{aligned}
	\right.
\end{equation}
where $a_{pi}$, $b_{i}$ and $c_{i}$ are coefficients corresponding to the explicit $p$th-order RK method.

We use RFM to represent the approximate solution $\tilde{\boldsymbol{\phi}}_{t_{k+1}}(\boldsymbol{x})$ of $\boldsymbol{\phi}_{t_{k+1}}(\boldsymbol{x})$
\begin{equation}
\begin{aligned}
\tilde{\boldsymbol{\phi}}_{t_{k+1}}(\boldsymbol{x})=&\left(\sum_{n=1}^{M_p} \psi_n(\boldsymbol{x}) \sum_{j=1}^{J_n} u_{n j, t_{k+1}}^1 \varphi_{n j, t_{k+1}}^1(\boldsymbol{x}),\cdots,\right.\\
&\left.\sum_{n=1}^{M_p} \psi_n(\boldsymbol{x}) \sum_{j=1}^{J_n} u_{n j, t_{k+1}}^{d_{\phi}} \varphi_{n j, t_{k+1}}^{d_{\phi}}(\boldsymbol{x})\right)^T.
\end{aligned}\label{eq:multiphaseapproximatesolution}
\end{equation}
% where $M_p=N_yN_x$ is the number of sub-domains of $\Omega$ and $J_n$ is the number of basis functions in each sub-domain.\\
Similarly, we sample $Q$ collocation points $\left\{\boldsymbol{x}_{q}^{n}\right\}$ in each $\Omega_n$ and compute the rescaling parameters
$\boldsymbol{\lambda}_{n, q, t_{k+1}}=diag(\lambda_{n, q, t_{k+1}}^{1},\cdots,\lambda_{n, q, t_{k+1}}^{d_{\phi}})$ and $\boldsymbol{\lambda}_{n, b, t_{k+1}}=diag(\lambda_{n, b, t_{k+1}}^{1},\cdots,\lambda_{n, b, t_{k+1}}^{d_{\phi}})$ as in \eqref{eq:RFMlambda}. Then, the RK-RFM minimizing the loss

\begin{equation}
	\begin{aligned}
&\operatorname{Loss}\left(\left\{u_{n j, t_{k+1}}^{i}\right\}\right)\\
=&\sum_{n=1}^{M_p}\left(\sum_{q=1}^Q\left\|\boldsymbol{\lambda}_{n, q, t_{k+1}}\left(\tilde{\boldsymbol{\phi}}_{t_{k+1}}\left(\boldsymbol{x}_{q}^{n}\right)-\boldsymbol{\phi}_{t_{k+1}}\left(\boldsymbol{x}_{q}^{n}\right)\right)\right\|_2^2\right)\\
&+\sum_{n=1}^{M_p}\left(\sum_{\boldsymbol{x}_{q}^{n} \in \partial \Omega}\left\|\boldsymbol{\lambda}_{n, b, t_{k+1}}\left(\mathcal{B}\tilde{\boldsymbol{\phi}}_{t_{k+1}}\left(\boldsymbol{x}_{q}^{n}\right)-\boldsymbol{g}\left(\boldsymbol{x}_{q}^{n},t_{k+1}\right)\right)\right\|_2^2\right)
\end{aligned}\label{eq:multiphaseloss}
\end{equation}
to find the optimal coefficients $\boldsymbol{U}_{t_{k+1}}=\left(u_{nj, t_{k+1}}^i\right)^{\top}$ of $\tilde{\boldsymbol{\phi}}_{t_{k+1}}(\boldsymbol{x})$ defined in \eqref{eq:multiphaseapproximatesolution}. In this article, we choose to use the PoU function $\psi$ defined in \eqref{eq:psi1}. Therefore, regularization terms should be added into the loss function \eqref{eq:multiphaseloss} to ensure continuity between subdomains. Specifically, following the approach in \cite{dong2021local}, we impose $C^1$
continuity conditions at the interfaces between subdomains. The loss functional \eqref{eq:multiphaseloss} implies 
\begin{equation}
\mathcal{A}_{t_{k+1}}\boldsymbol{U}_{t_{k+1}}=\boldsymbol{f}_{t_{k+1}},
\label{eq:multiphaseloss1}
\end{equation}
where $\mathcal{A}_{t_{k+1}}$ is the coefficient matrix related to $\tilde{\boldsymbol{\phi}}_{t_{k+1}}\left(\boldsymbol{x}_{q}^{n}\right)$ and  $\boldsymbol{f}_{t_{k+1}}$ is the right-hand side term related to $\boldsymbol{\phi}_{t_{k+1}}\left(\boldsymbol{x}_{q}^{n}\right)$ and $\boldsymbol{g}\left(\boldsymbol{x}_{q}^{n},t_{k+1}\right)$ (see Appendix \ref{app:A}). Thus, we can obtain the the optimal $\boldsymbol{U}_{t_{k+1}}$ by solving \eqref{eq:multiphaseloss1} using the linear least-squares method. Once the optimal $\boldsymbol{U}_{t_{k+1}}$ is obtained, the
$\tilde{\boldsymbol{\phi}}_{t_{k+1}}(\boldsymbol{x})$  is used as the initial condition for the next time step in the RK method. We summarize the steps of the RK-RFM in Algorithm \ref{alg:Algorithm 1} and depict the architecture in  Figure \ref{fig:fig1}.

\begin{remark}(Without automatic differentiation)
It should be noted that, with the PoU function $\psi$ defined in \eqref{eq:psi1}, the calculation of the derivatives 
\begin{equation}
	\begin{aligned}
&\nabla_{\boldsymbol{x}}^{l} \tilde{\boldsymbol{\phi}}_{t_{k+1}}(\boldsymbol{x})\\
= &\left(\sum_{n=1}^{M_p} \psi_n(\boldsymbol{x}) \sum_{j=1}^{J_n} u_{n j, t_{k+1}}^1 \nabla_{\boldsymbol{x}}^{l}\varphi_{n j, t_{k+1}}^1(\boldsymbol{x}),\cdots,\right.\\
&\left.\sum_{n=1}^{M_p} \psi_n(\boldsymbol{x}) \sum_{j=1}^{J_n} u_{n j, t_{k+1}}^{d_{\phi}} \nabla_{\boldsymbol{x}}^{l}\varphi_{n j, t_{k+1}}^{d_{\phi}}(\boldsymbol{x})\right)^T, \ l=1,\cdots,mp,
	\end{aligned}\label{eq:derivatives}
 \end{equation}
is essential when dealing with the $m$th-order equation  in the explicit $p$th-order  RK method. When the number of subdomains $M_p$, the number of random feature functions $J_n$, the number of sample points $Q$,  the order of the RK method $p$ or the order of the  equation $m$ are large, using automatic differentiation can consume a large amount of computational resources and take a long time when we calculate \eqref{eq:derivatives}.  It is better to manually derive the necessary derivative functions
\begin{equation}
	\begin{aligned}
\nabla_{\boldsymbol{x}}^{l}\varphi(\boldsymbol{x},\boldsymbol{W},b)=\nabla_{\boldsymbol{x}}^{l}\sigma\left(\boldsymbol{W}\cdot \tilde{\boldsymbol{x}}+b\right), \quad l=1,\cdots,mp,
\end{aligned}\label{eq:derivativefunction}
 \end{equation}
 before starting the computation and store them. By this way, we can directly obtain $\nabla_{\boldsymbol{x}}^{l}\varphi_{n j, t_{k+1}}^{i}(\boldsymbol{x}) = \nabla_{\boldsymbol{x}}^{l}\varphi(\boldsymbol{x},\boldsymbol{W}_{n j, t_{k+1}}^{i},b_{n j, t_{k+1}}^{i})$ by using \eqref{eq:derivativefunction}  when needed, which will significantly reduce the use of computational resources and improve computation efficiency.
\end{remark}

\begin{remark}(Parallelization)
In the aforementioned process, we need to generate different random feature functions $\varphi^i_{nj, t_{k+1}}$ for different $\tilde{\boldsymbol{\phi}}_{t_{k+1}}^i$ and calculate the corresponding coefficient matrices and right-hand side terms. However, due to the limitations of computational resources, we need to solve for $\boldsymbol{U}_{t_{k+1}}^i$ using the linear least squares method sequentially. When the dimension $d_{\phi}$ of $\tilde{\boldsymbol{\phi}}_{t_{k+1}}$ is very large, this process becomes very time-consuming. To improve the computational efficiency and reduce the consumption of computational resources, we choose to use the same random feature functions $\varphi_{nj, t_{k+1}}$ for different $\tilde{\boldsymbol{\phi}}_{t_{k+1}}^i$. This means that in each time step $t_{k+1}$, we only need to calculate the coefficient matrix $ \mathcal{A}_{t_{k+1}}$ once, and the right-hand 
side terms $\boldsymbol{f}^i_{t_{k+1}}$ corresponding to different $\phi_i$ at time $t_{k+1}$ are assembled into a matrix. Specifically, we have  
\begin{equation}
\begin{aligned}
\mathcal{A}_{t_{k+1}}\boldsymbol{U}_{t_{k+1}}&=\mathcal{A}_{t_{k+1}}\left[\boldsymbol{U}^1_{t_{k+1}},\cdots,\boldsymbol{U}^{d_{\phi}}_{t_{k+1}}\right]\\
&=\left[\boldsymbol{f}^1_{t_{k+1}},\cdots,\boldsymbol{f}^{d_{\phi}}_{t_{k+1}}\right]\\
&=\boldsymbol{f}_{t_{k+1}}.
\end{aligned}\label{eq:matrix_parallelization}
\end{equation}
Then, we only need to use the linear least squares method once to obtain the optimal  coefficients $\boldsymbol{U}_{t_{k+1}}$  corresponding to $\tilde{\boldsymbol{\phi}}_{t_{k+1}}$.
\end{remark}

\begin{figure}[htp]
\begin{algorithm}[H]
	\caption{ The Runge-Kutta Random Feature Method}
	\label{alg:Algorithm 1}
	\begin{algorithmic}[1]
		\Require
		 number of the subdomains $M_p$; number of the random feature functions on each subdomain $J_n$;  number of collocation points in each subdomain $Q$; terminal time $T$;
        number of time intervals $K$;  coefficients in the explicit $p$th-order RK method $a_{pi}$, $b_{i}$ and $c_{i}$;  range of uniform distribution $R_m$; constant coefficient of the rescaling parameters $c$.
		\Ensure
		The optimal $\boldsymbol{U}_{t_1},\cdots,\boldsymbol{U}_{t_K}$;
		\State Divide $\Omega$ into $M_p$ non-overlapping subdomains ${\Omega_n}$;
        \State  Sample $Q$ collocation points $\left\{ \boldsymbol{x}_{q}^{n}\right\}$ in each $\Omega_n$;
        \State  Compute the derivative functions \eqref{eq:derivativefunction}; 
		\State  \textbf{for} $k=0,1,2, \cdots, K-1$ \textbf{do}
		\State \quad Construct $J_n $ random feature functions $\varphi_{n j, t_{k+1}}$ on $\Omega_n$ by the uniform distribution $\mathbb{U}\left(-R_m, R_m\right)$;
        \State  \quad Compute $\boldsymbol{\phi}_{t_{k+1}}(\boldsymbol{x})$ according to \eqref{eq:porderrk};
        \State  \quad Assemble  the loss matrix \eqref{eq:matrix_parallelization} and solve it by the linear least-squares method to obtain the optimal $\boldsymbol{U}_{t_{k+1}}$;
        \State  \quad Compute $\tilde{\boldsymbol{\phi}}_{t_{k+1}}(\boldsymbol{x})$ according to \eqref{eq:multiphaseapproximatesolution};
        \State  \quad $\boldsymbol{\phi}_{t_{k+1}}(\boldsymbol{x}):= \tilde{\boldsymbol{\phi}}_{t_{k+1}}(\boldsymbol{x})$
        \State \textbf{end}
        \State Return the optimal $\boldsymbol{U}_{t_1},\cdots,\boldsymbol{U}_{t_K}$.
	\end{algorithmic}
\end{algorithm}
\end{figure}

\begin{figure*}[htp]
	\center{\includegraphics[scale=0.26]  {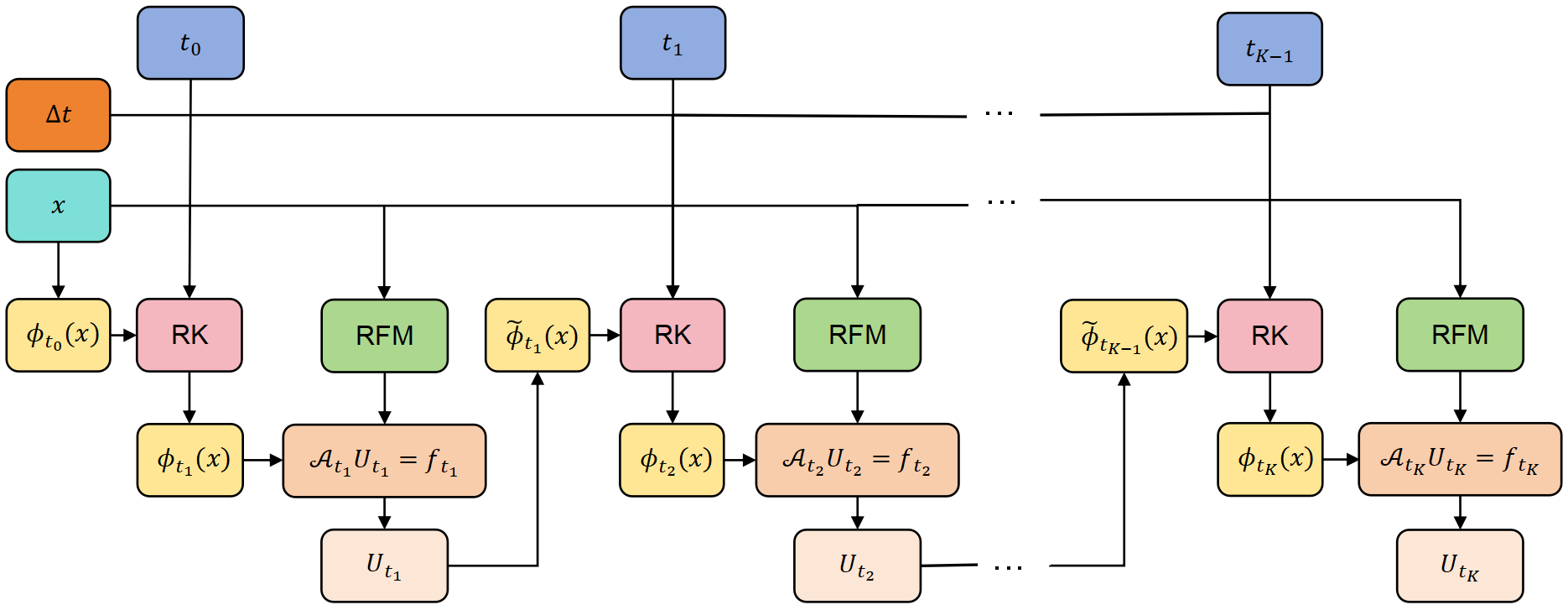}} 
	\caption{The architecture of of Algorithm \ref{alg:Algorithm 1}.}
	\label{fig:fig1}
\end{figure*}
\subsection{Error Estimates of  Runge-Kutta Random Feature Method}
\label{sec03:sec02}
\begin{lemma} \label{lemma1}
 Let $\sigma$ be a bounded nonconstant piecewise continuous activation function, $d_{\phi}=1$, $G_{M}(\boldsymbol{x})$ be of the form \eqref{eq:approximatesolution}, and $M = M_pJ_n$ denotes the degree of freedom. Then, for any $f \in C(\Omega)$, there exists a sequence of $G_{M}(\boldsymbol{x})$, such that
\begin{equation*}
    \begin{aligned}
% \lim_{J_n\to\infty}\left\|G_{M_p,J_n}(\boldsymbol{x})-f(\boldsymbol{x})\right\|_{L^2}=0, \quad
\lim_{M\to\infty}\left\|G_{M}(\boldsymbol{x})-f(\boldsymbol{x})\right\|_{L^2}=0, \quad \forall \boldsymbol{x} \in \Omega.
\end{aligned}
\end{equation*}
% if the optimal  coefficients $\boldsymbol{U}=\left(u_{nj}\right)^{\top}$  are determined by ordinary least square to minimize $\left\|G_{M_p,J_n}(\boldsymbol{x})-f(\boldsymbol{x})\right\|_{L^2}$
\begin{proof}
It is a direct consequence of Theorem 2.1 in reference \cite{huang2006universal}.
%in %\cite[Corollary 2.2]{chen2023random}
\end{proof}
\end{lemma}
% \begin{lemma} \label{lemma1}
%  Let $\sigma$ be a bounded nonconstant piecewise continuous activation function, $d_{\phi}=1$, and $G_{M_p,J_n}(\boldsymbol{x})$ be of the form \eqref{eq:approximatesolution}.
% %  Then finite sums of
% % \begin{equation*}
% %     \begin{aligned}
% % G_M(\boldsymbol{x})=\sum_{j=1}^{J_n} u_j \sigma\left(\boldsymbol{W}_j\cdot \boldsymbol{x}+b_j\right)
% % \end{aligned}
% % \end{equation*}
% % are dense in $C(\Omega)$. 
%  For any $f \in C(\Omega)$, there exists a sequence of $G_{M_p,J_n}(\boldsymbol{x})$, such that
% \begin{equation*}
%     \begin{aligned}
% \lim_{J_n\to\infty}\left\|G_{M_p,J_n}(\boldsymbol{x})-f(\boldsymbol{x})\right\|=0, \quad
% \lim_{M_p\to\infty}\left\|G_{M_p,J_n}(\boldsymbol{x})-f(\boldsymbol{x})\right\|=0, \quad \forall \boldsymbol{x} \in \Omega.
% \end{aligned}
% \end{equation*}
% \begin{proof}
% It is a direct consequence of Corollary 2.2 in \cite{chen2023random}.
% %in %\cite[Corollary 2.2]{chen2023random}
% \end{proof}
% \end{lemma}

%In practice, however, 
%The numerical solution $\hat{\boldsymbol{U}}$ is denoted by $\hat{\boldsymbol{U}}=\left(\hat{u}_{j}\right)^{\top}$  and
%the optimal solution $\tilde{\boldsymbol{U}} $is denoted by $\tilde{\boldsymbol{U}}=\left(\tilde{u}_{j}\right)^{\top}$, and we denote the difference by $\delta u$. 
We use $\delta u$ to denote the difference between the numerical solution $\hat{\boldsymbol{U}} = \left(\hat{u}_{j}\right)^{\top}$ and the optimal solution $\tilde{\boldsymbol{U}} = \left(\tilde{u}_{j}\right)^{\top}$ in the subdomain $\Omega_{n}$.
Let $\left\{\boldsymbol{b}_j\right\}_{j=1}^{J_n}$ be a unitary orthogonal basis in $\mathbb{R}^{J_n}$, %Since $\left\{\boldsymbol{b}_j\right\}_{j=1}^{J_n}$ forms a basis, 
there exists $\left\{\delta u_j\right\}_{j=1}^{J_n}$ such that $\delta u=\sum_{j=1}^{J_n} \delta u_j \boldsymbol{b}_j$, where $\delta u_j, j=1, \cdots, J_n$ are independent and identically distributed random variables with 
\begin{eqnarray}\label{eq:definition1}
\mathbb{E}\left[\delta u_j\right]=0, \  \mathbb{E}\left[\vert\delta u_j\vert\right]=\mu>0, \ \mathbb{E}\left[\left(\delta u_j\right)^2\right]=\delta^2>0. 
\end{eqnarray}

\begin{lemma} \label{lemma2}
For the RFM, we denote $\hat{\boldsymbol{U}}=\left(\hat{u}_{nj}^i\right)^{\top}$ and $\tilde{\boldsymbol{U}}=\left(\tilde{u}_{nj}^i\right)^{\top}$ as the numerical solution and the optimal solution obtained by minimizing the loss function \eqref{eq:loss}, respectively. 
% Let $\delta u$ be the error between the numerical solution and the optimal solution to the least-squares problem.
 We have
\begin{equation*}
    \begin{aligned}
\mathbb{E}_{\delta u}\left[\left\|\hat{\boldsymbol{U}}-\tilde{\boldsymbol{U}}\right\|_{L^2}\right] \leq \sqrt{d_{\phi}M_p J_n } \delta,
\end{aligned}
\end{equation*}
where $\delta$ is defined in \eqref{eq:definition1}.
\begin{proof}
\begin{equation*}
\begin{aligned}
\left(\mathbb{E}_{\delta u}\left[\left\|\hat{\boldsymbol{U}}-\tilde{\boldsymbol{U}}\right\|_{L^2}\right]\right)^2 & \leq \mathbb{E}_{\delta u}\left[\left\|\hat{\boldsymbol{U}}-\tilde{\boldsymbol{U}}\right\|^2_{L^2}\right]\\
& =\mathbb{E}_{\delta u}\left[\sum_{i=1}^{d_{\phi}}\sum_{n=1}^{M_p} \left\|\sum_{j=1}^{J_n} \delta u_{nj}^i \boldsymbol{b}_j\right\|^2_{L^2}\right] \\
&=d_{\phi}M_p J_n \delta^2.
\end{aligned}
\end{equation*}
\end{proof}
\end{lemma}

\begin{lemma} \label{lemma3} 
Let two functions $\hat{\boldsymbol{G}}_M(\boldsymbol{x})$ and $\tilde{\boldsymbol{G}}_M(\boldsymbol{x})$ be of the formulation \eqref{eq:approximatesolution}  and the coefficients being $\hat{\boldsymbol{U}}=\left(\hat{u}_{nj}^i\right)^{\top}, \tilde{\boldsymbol{U}}=\left(\tilde{u}_{nj}^i\right)^{\top}$, respectively, then there exists $C_1 \geq C_2>0$ such that
\begin{equation*}
\begin{aligned}
C_1\|\hat{\boldsymbol{U}}-\tilde{\boldsymbol{U}}\|_{L^2} \geq\left\|\hat{\boldsymbol{G}}_M(\boldsymbol{x})-\tilde{\boldsymbol{G}}_M(\boldsymbol{x})\right\|_{L^2} \geq C_2\|\hat{\boldsymbol{U}}-\tilde{\boldsymbol{U}}\|_{L^2}.
\end{aligned}
\end{equation*}

\begin{proof}
Notice that
\begin{equation*}
\begin{aligned}
&\left\|\hat{\boldsymbol{G}}_M(\boldsymbol{x})-\tilde{\boldsymbol{G}}_M(\boldsymbol{x})\right\|_{L^2}\\
&=\left(\sum_{i=1}^{d_{\phi}}\sum_{n=1}^{M_p} \sum_{j=1}^{J_n}\left(\hat{u}_{nj}^i-\tilde{u}_{nj}^i\right)^2 \int_{\Omega} \varphi_{nj}^i(\boldsymbol{x})^2 \mathrm{d}\boldsymbol{x}\right)^{\frac{1}{2}} \\
&\leq C_1\|\hat{\boldsymbol{U}}-\tilde{\boldsymbol{U}}\|_{L^2},
\end{aligned}
\end{equation*}
where $C_1=\max _{n, j, i}\left\|\varphi_{nj}^i(\boldsymbol{x})\right\|_{L^2}$, and
\begin{equation*}
\begin{aligned}
&\left\|\hat{\boldsymbol{G}}_M(\boldsymbol{x})-\tilde{\boldsymbol{G}}_M(\boldsymbol{x})\right\|_{L^2}\\
&=\left(\sum_{i=1}^{d_{\phi}}\sum_{n=1}^{M_p} \sum_{j=1}^{J_n}\left(\hat{u}_{nj}^i-\tilde{u}_{nj}^i\right)^2 \int_{\Omega } \varphi_{nj}^i(\boldsymbol{x})^2  \mathrm{d}\boldsymbol{x} \right)^{\frac{1}{2}} \\
&\geq C_2\|\hat{\boldsymbol{U}}-\tilde{\boldsymbol{U}}\|_{L^2},
\end{aligned}
\end{equation*}
where $C_2=\min _{n, j, i}\left\|\varphi_{nj}^i(\boldsymbol{x})\right\|_{L^2}$. Moreover, for all $\varphi^i_{nj} \in \mathcal{L}^2(\Omega)$, we have $C_1, C_2<+\infty$.
\end{proof}
\end{lemma}

\begin{thm}\label{thm1}
Let $\boldsymbol{\phi}^{e}_{t_{k}}$ be the exact solution of the equation and $\boldsymbol{\phi}_{t_{k}}$ be the solution to  explicit $p$th-order RK method \eqref{eq:porderrk}. Let $\sigma$ be a bounded nonconstant piecewise continuous activation function,  $\tilde{\boldsymbol{\phi}}_{t_{k}}$ and $\hat{\boldsymbol{\phi}}_{t_{k}}$ in the form of \eqref{eq:multiphaseapproximatesolution} are the optimal solution and  the numerical solution obtained by minimizing the least-squares problem \eqref{eq:matrix_parallelization} ,respectively. 
% Let $\delta u$ be the error between the numerical solution and the optimal solution to the least-squares problem. 
We assume there exists $L>0$ such that
\begin{equation*}
\begin{aligned}
&\left\|H(\boldsymbol{x}, t_k, \hat{\boldsymbol{\phi}}_{t_{k}}(\boldsymbol{x}), \Delta t)-H(\boldsymbol{x}, t_k, \boldsymbol{\phi}^{e}_{t_{k}}(\boldsymbol{x}), \Delta t)\right\|_{L^2}\\ 
&\leq L \left\|\hat{\boldsymbol{\phi}}_{t_{k}}(\boldsymbol{x})-\boldsymbol{\phi}^{e}_{t_{k}}(\boldsymbol{x})\right\|_{L^2}, \\
&\forall \hat{\boldsymbol{\phi}}_{t_{k}}(\boldsymbol{x}), \boldsymbol{\phi}^{e}_{t_{k}}(\boldsymbol{x})\in R^{d_{\phi}}, \boldsymbol{x}\in\Omega, k=0,\cdots,K, 
\end{aligned}
\end{equation*}
where $H$ is defined in \eqref{eq:porderrk} and given $\epsilon>0$,  we also assume that
\begin{equation*}
\begin{aligned}
&\left\|\tilde{\boldsymbol{\phi}}_{t_{k}}(\boldsymbol{x})-\boldsymbol{\phi}_{t_{k}}(\boldsymbol{x})\right\|_{L^2} \leq\sqrt{d_{\phi}} \epsilon, \\
&\forall \boldsymbol{x}\in\Omega, k=0,\cdots,K, n=1, \cdots, M_p.  
\end{aligned}
\end{equation*}
Then, there exists  $C_1, C_2>0$ such that
\begin{equation}
\begin{aligned}
&\mathbb{E}_{\delta u}\left[\left\|\hat{\boldsymbol{\phi}}_{t_{k}}(\boldsymbol{x})-\boldsymbol{\phi}^{e}_{t_{k}}(\boldsymbol{x})\right\|_{L^2}\right] \\
&\leq \frac{e^{L(t_{k}-t_{0})}-1}{\Delta t L}(C_1\sqrt{d_{\phi}M_p J_n } \delta+\sqrt{d_{\phi}} \epsilon)+C_2\Delta t^p \\
&= e_1(\Delta t)+e_2(\Delta t),\  \forall \boldsymbol{x}\in\Omega, k=0,\cdots,K,
\end{aligned}\label{eq:error}
\end{equation}
where $e_1(\Delta t) = \frac{e^{L(t_{k}-t_{0})}-1}{\Delta t L}(C_1\sqrt{d_{\phi}M_p J_n } \delta+\sqrt{d_{\phi}} \epsilon)$ and $e_2(\Delta t) = C_2\Delta t^p$.
\begin{proof} Notice that
\begin{equation*}
\begin{aligned}
&\mathbb{E}_{\delta u}\left[\left\|\hat{\boldsymbol{\phi}}_{t_{k}}(\boldsymbol{x})-\boldsymbol{\phi}^{e}_{t_{k}}(\boldsymbol{x})\right\|_{L^2}\right] \\
\leq& \mathbb{E}_{\delta u}\left[\left\|\hat{\boldsymbol{\phi}}_{t_{k}}(\boldsymbol{x})-\tilde{\boldsymbol{\phi}}_{t_{k}}(\boldsymbol{x})\right\|_{L^2}\right]+\left\|\tilde{\boldsymbol{\phi}}_{t_{k}}(\boldsymbol{x})-\boldsymbol{\phi}_{t_{k}}(\boldsymbol{x})\right\|_{L_2}\\
&+\mathbb{E}_{\delta u}\left[\left\|\boldsymbol{\phi}_{t_{k}}(\boldsymbol{x})-\boldsymbol{\phi}^{e}_{t_{k}}(\boldsymbol{x})\right\|_{L_2}\right].
% & \leq \alpha \sigma \sqrt{N_t}+\left(\sum_{n=1}^{N_t} \epsilon^2\right)^{\frac{1}{2}} \\
% & =\alpha \sigma \sqrt{N_t}+\sqrt{N_t} \epsilon.
\end{aligned}
\end{equation*}
Therefore, we have
\begin{equation*}
\begin{aligned}
&\left\|\boldsymbol{\phi}_{t_{k}}(\boldsymbol{x})-\boldsymbol{\phi}^{e}_{t_{k}}(\boldsymbol{x})\right\|_{L_2}\\
=&\left\|\boldsymbol{\phi}_{t_{k}}(\boldsymbol{x})-(\boldsymbol{\phi}^{e}_{t_{k-1}}(\boldsymbol{x})+\Delta t H(\boldsymbol{x},t_{k-1},\boldsymbol{\phi}^{e}_{t_{k-1}}(\boldsymbol{x}),\Delta t))\right.\\
&\left.+(\boldsymbol{\phi}^{e}_{t_{k-1}}(\boldsymbol{x})+\Delta t H(\boldsymbol{x},t_{k-1},\boldsymbol{\phi}^{e}_{t_{k-1}}(\boldsymbol{x}),\Delta t))-\boldsymbol{\phi}^{e}_{t_{k}}(\boldsymbol{x})\right\|_{L_2}\\
=& \left\|(\hat{\boldsymbol{\phi}}_{t_{k-1}}(\boldsymbol{x})+\Delta t H(\boldsymbol{x},t_{k-1},\hat{\boldsymbol{\phi}}_{t_{k-1}}(\boldsymbol{x}),\Delta t))\right.\\
&\left.-(\boldsymbol{\phi}^{e}_{t_{k-1}}(\boldsymbol{x})+\Delta t H(\boldsymbol{x},t_{k-1},\boldsymbol{\phi}^{e}_{t_{k-1}}(\boldsymbol{x}),\Delta t))\right.\\
&\left.+(\boldsymbol{\phi}^{e}_{t_{k-1}}(\boldsymbol{x})+\Delta t H(\boldsymbol{x},t_{k-1},\boldsymbol{\phi}^{e}_{t_{k-1}}(\boldsymbol{x}),\Delta t))-\boldsymbol{\phi}^{e}_{t_{k}}(\boldsymbol{x})\right\|_{L_2}\\
\leq&\left\|\hat{\boldsymbol{\phi}}_{t_{k-1}}(\boldsymbol{x})-\boldsymbol{\phi}^{e}_{t_{k-1}}(\boldsymbol{x})\right\|_{L_2}\\
&+ \Delta t  \left\|H(\boldsymbol{x},t_{k-1},\hat{\boldsymbol{\phi}}_{t_{k-1}}(\boldsymbol{x}),\Delta t)- H(\boldsymbol{x},t_{k-1},\boldsymbol{\phi}^{e}_{t_{k-1}}(\boldsymbol{x}),\Delta t)\right\|_{L_2}\\
& +\left\|(\boldsymbol{\phi}^{e}_{t_{k-1}}(\boldsymbol{x})+\Delta t H(\boldsymbol{x},t_{k-1},\boldsymbol{\phi}^{e}_{t_{k-1}}(\boldsymbol{x}),\Delta t))-\boldsymbol{\phi}^{e}_{t_{k}}(\boldsymbol{x})\right\|_{L_2}\\
\leq & (1+\Delta t L)\left\|\hat{\boldsymbol{\phi}}_{t_{k-1}}(\boldsymbol{x})-\boldsymbol{\phi}^{e}_{t_{k-1}}(\boldsymbol{x})\right\|_{L_2}+C\Delta t^{p+1}.
\end{aligned}
\end{equation*}
Then
\begin{equation*}
\begin{aligned}
&\mathbb{E}_{\delta u}\left[\left\|\hat{\boldsymbol{\phi}}_{t_{k}}(\boldsymbol{x})-\boldsymbol{\phi}^{e}_{t_{k}}(\boldsymbol{x})\right\|_{L^2}\right] \\
\leq& C_1\sqrt{d_{\phi}M_p J_n } \delta + \sqrt{d_{\phi}} \epsilon+C_2\Delta t^{p+1} \\
&+ (1+\Delta t L)\mathbb{E}_{\delta u}\left[\left\|\hat{\boldsymbol{\phi}}_{t_{k-1}}(\boldsymbol{x})-\boldsymbol{\phi}^{e}_{t_{k-1}}(\boldsymbol{x})\right\|_{L^2}\right]\\
=&(1+(1+\Delta t L)+\cdots+(1+\Delta t L)^{k-1})(C_1\sqrt{d_{\phi}M_p J_n } \delta\\
&+\sqrt{d_{\phi}} \epsilon+C_2\Delta t^{p+1})\\
&+(1+\Delta t L)^{k}\mathbb{E}_{\delta u}\left[\left\|\hat{\boldsymbol{\phi}}_{t_{0}}(\boldsymbol{x})-\boldsymbol{\phi}^{e}_{t_{0}}(\boldsymbol{x})\right\|_{L^2}\right]\\
=&\frac{(1+\Delta t L)^{k}-1}{\Delta t L}(C_1\sqrt{d_{\phi}M_p J_n } \delta+\sqrt{d_{\phi}} \epsilon+C_2\Delta t^{p+1})\\
\leq& \frac{(e^{\Delta t L})^k-1}{\Delta t L}(C_1\sqrt{d_{\phi}M_p J_n } \delta+\sqrt{d_{\phi}} \epsilon+C_2\Delta t^{p+1})\\
=&\frac{e^{L(t_{k}-t_{0})}-1}{\Delta t L}(C_1\sqrt{d_{\phi}M_p J_n } \delta+\sqrt{d_{\phi}} \epsilon)+C_2\Delta t^p.
% & \leq \alpha \sigma \sqrt{N_t}+\left(\sum_{n=1}^{N_t} \epsilon^2\right)^{\frac{1}{2}} \\
% & =\alpha \sigma \sqrt{N_t}+\sqrt{N_t} \epsilon.
\end{aligned}
\end{equation*}
\end{proof}
\end{thm}

\begin{remark}(The choice of $\Delta t$)
Note that the error \eqref{eq:error} consists of two parts, %$e_1(\Delta t) $ and $e_2(\Delta t)$,
where $e_1(\Delta t)$ represents the error caused by RFM and optimization error, which is inversely proportional to $\Delta t$, the $e_2(\Delta t)$ represents the error introduced by the explicit $p$th-order RK method \eqref{eq:porderrk}, which is directly proportional to $\Delta t^p$. A larger $\Delta t$ can lead to a larger $e_2(\Delta t)$ and smaller $e_1(\Delta t)$, while a smaller $\Delta t$ may result in a smaller $e_2(\Delta t)$ and larger $e_1(\Delta t)$, also costs much computational time.
%Conversely, choosing a smaller $\Delta t$ may not necessarily have advantages. Although a smaller $\Delta t$ results in a smaller $e_2(\Delta t)$, it increases $e_1(\Delta t)$ and adds to the computational time cost. 
Ideally, $\Delta t$ should be selected to satisfy $e_1(\Delta t) \approx  e_2(\Delta t)$. However, this principle may not always be feasible since the values of $\delta$ and $\epsilon$ are often unknown. Regardless, this principle still provides valuable guidance for the selection of the time step $\Delta t$.
\end{remark}

% \begin{remark}(The choice of $M_p$, $J_n$ and $Q$)
% Note that the error $e_1(\Delta t)$ consists of two parts, 
%  the error caused by RFM $\frac{e^{L(t_{k}-t_{0})}-1}{\Delta t L}\sqrt{d_{\phi}}\epsilon$ and the optimization error $\frac{e^{L(t_{k}-t_{0})}-1}{\Delta t L}(C_1\sqrt{d_{\phi}M_p J_n } \delta)$. From the Lemma \ref{lemma1}, it is known that increasing $M_p$ and $J_n$ can reduce the error caused by RFM, but according to the Theorem \ref{thm1}, the optimization error will correspondingly increase. Typically, increasing $Q$ can cause $\delta$ to decrease, so choosing larger $M_p$, $J_n$ and $Q$ can yield better results.
% % Ideally, $M_p$, $J_n$ and $Q$ should be selected to satisfy $\epsilon \approx  C_1\sqrt{M_p J_n } \delta$. 
% % However, this principle may not always be feasible since the values of $\delta$ and $\epsilon$ are often unknown. 
% % Typically, the magnitude of $\delta$ is  smaller than that of $\epsilon$, and choosing larger $M_p$ and $J_n$ can yield better results.
% \end{remark}

\section{Numerical experiments}
\label{sec04}
\begin{figure*}[htp]
	\center
        {\subfigure[The arrangement of subdomains. ]{		 
            \includegraphics[scale=0.18]{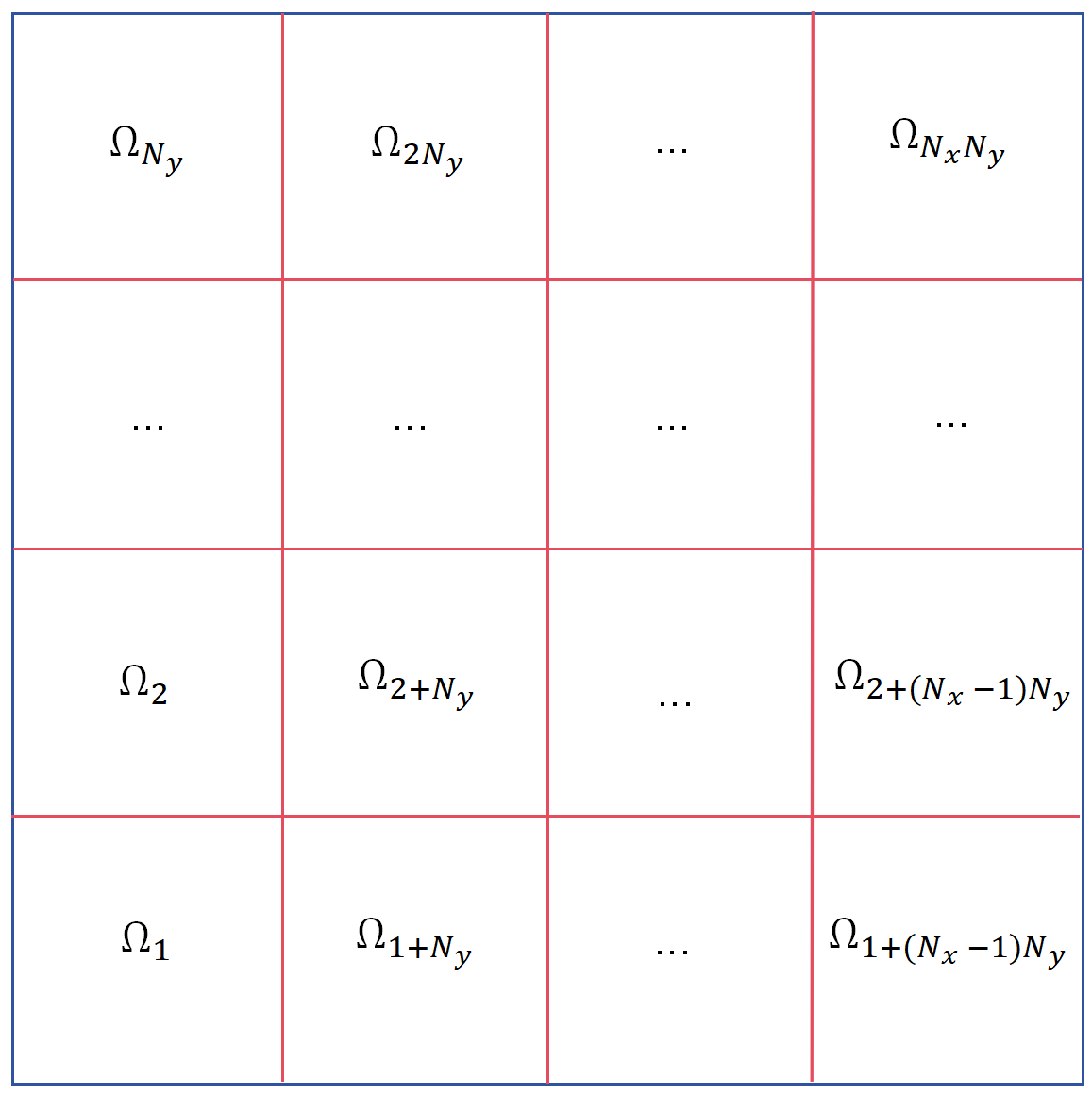}}
            \subfigure[The distribution of collocation points in subdomain $\Omega_{n}$. ]{		 
            \includegraphics[scale=0.19]{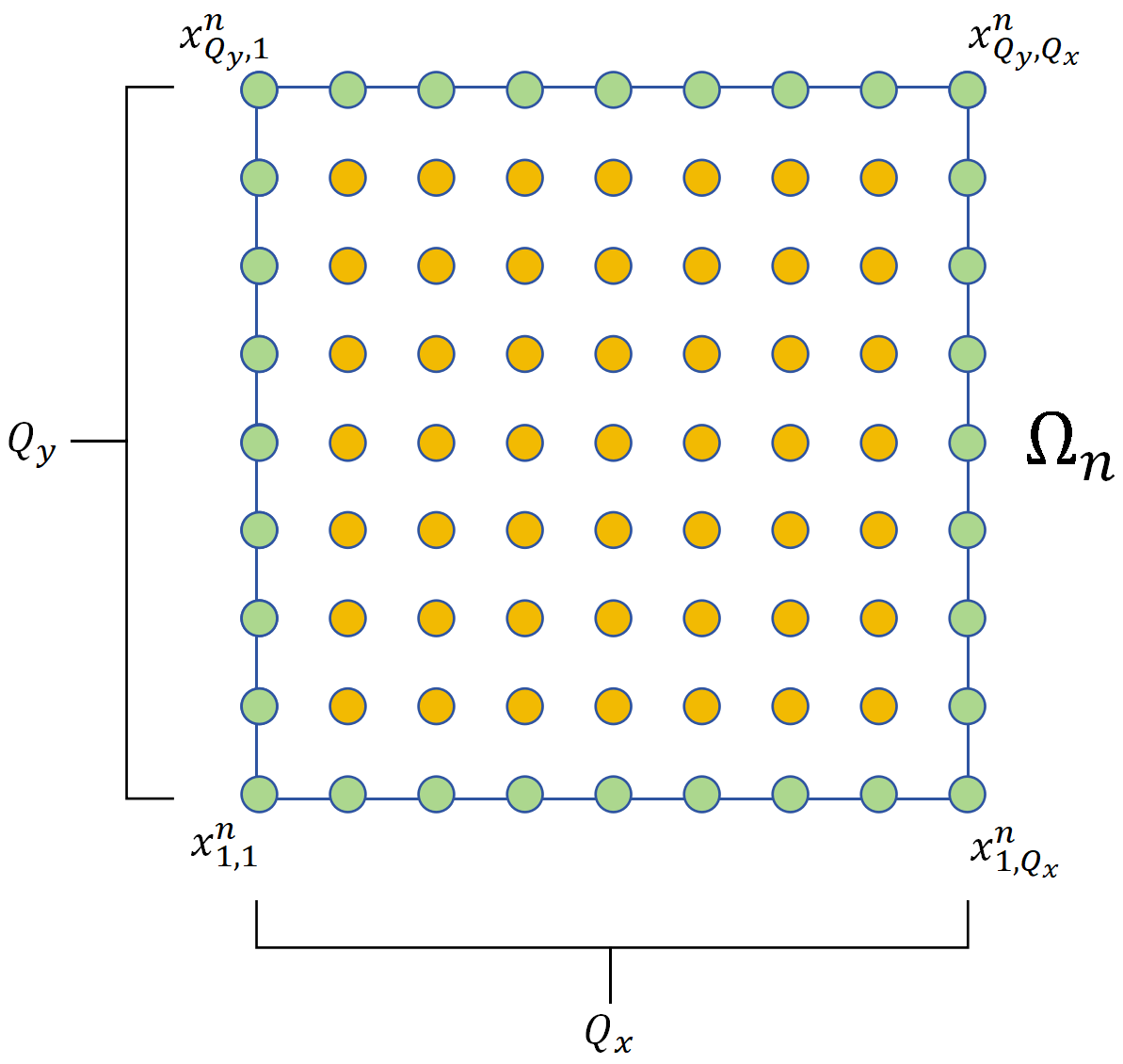}}}
	\caption{The arrangement of subdomains and the distribution of collocation points in subdomain $\Omega_{n}$. (a) The blue lines and red lines represent the boundaries and the interfaces between subdomains, respectively. (b) The green points and yellow points represent the boundary points and interior points in each subdomain, respectively.}
	\label{fig:fig13}
\end{figure*}

In this section, we present some numerical experiments to validate our methods. For quantitative comparison, we compare the errors of the numerical solution $\hat{\boldsymbol{\phi}}_{t_{K}}(\boldsymbol{x})$ and the exact solution $\boldsymbol{\phi}^{e}_{t_{K}}(\boldsymbol{x})$ in the relative $L^{\infty}$ norm and relative $L^2$ norm, which are defined as
\begin{equation*}
\begin{aligned}
&\vert\vert e\vert\vert_{L^{\infty}} = \frac{\vert\vert \hat{\boldsymbol{\phi}}_{t_{K}}(\boldsymbol{x})-\boldsymbol{\phi}^{e}_{t_{K}}(\boldsymbol{x})\vert\vert_{L^{\infty}}}{\vert\vert \boldsymbol{\phi}^{e}_{t_{K}}(\boldsymbol{x})\vert\vert_{L^{\infty}}},\\
 &\vert\vert e\vert\vert_{L^2}=\frac{\vert\vert \hat{\boldsymbol{\phi}}_{t_{K}}(\boldsymbol{x})-\boldsymbol{\phi}^{e}_{t_{K}}(\boldsymbol{x})\vert\vert_{L^{2}}}{\vert\vert \boldsymbol{\phi}^{e}_{t_{K}}(\boldsymbol{x})\vert\vert_{L^{2}}}.
 \end{aligned}
\end{equation*}
In two dimension, we set $M_p = N_xN_y$ and $Q=Q_xQ_y$,
where $N_x$ and $N_y$ represent the number of subdomains along the $x$
and $y$ directions, respectively, and $Q_x$ and $Q_y$ are the number of points in the $x$ and $y$ directions of the  subdomain $\Omega_{n}$. The arrangement of subdomains and the distribution of collocation points in subdomain $\Omega_{n}$ are shown in Figure \ref{fig:fig13}.

% \begin{figure}[htp]
% 	\centering 	 
% 		\includegraphics[scale=0.2]{Q_01.png}
% 	\caption{The distribution of collocation points on subdomain $\Omega_{n}$. The green points and yellow points represent the boundary points and interior points on each subdomain, respectively.}
% 	\label{fig:fig14}
% \end{figure}
 In the linear least squares simulations, we use the least squares solver \textit{torch.linalg.lstsq} available in PyTorch.
% We use the  standard Kaiming uniform initializer \textit{torch.nn.init.kaiming\_uniform\_} in Pytorch to generate $\left\{\boldsymbol{W}_{nj}\right\}$ and $\left\{b_{nj}\right\}$, where the parameter settings are \textit{a=1}, \textit{mode='fan\_in'} and \textit{nonlinearity='leaky\_relu'}.
% We use the  standard least squares solver \textit{torch.linalg.lstsq} in Pytorch to solve all the least squares problems. 
 We set the constant coefficient of the rescaling parameters $c = 100$. The second order explicit RK method ($p=2$) is utilized. In our experiments, we use AMD EPYC 7T83 CPU and NVIDIA RTX 4090 GPU for simulation. 
 \subsection{Convergence Test}
\label{sec06:sec01}
\begin{figure}[htp]
	\centering 	 
		\includegraphics[scale=0.34]{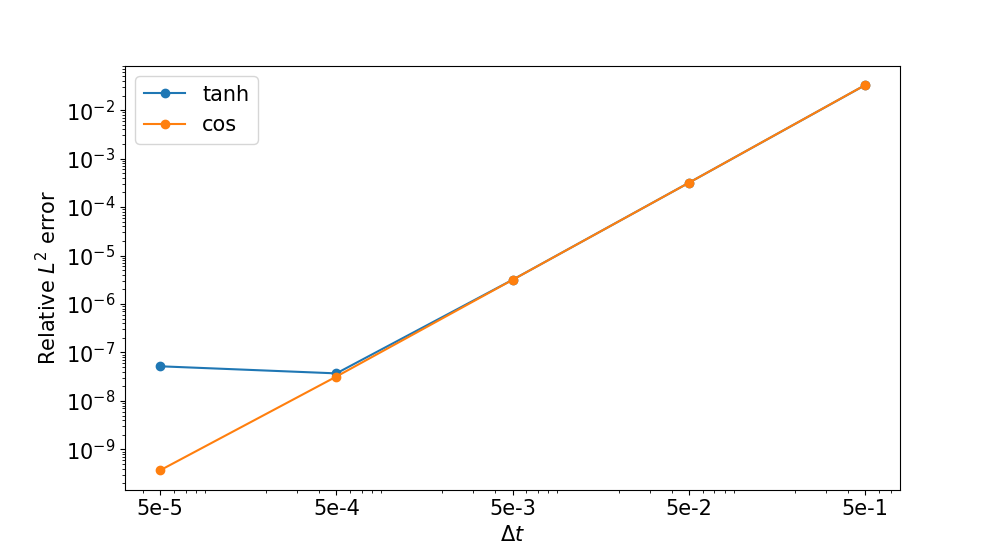}
	\caption{Relative $L^2$ error for \eqref{eq:exactsolution} obtained by the Algorithm \ref{alg:Algorithm 1} with $\Delta t=$ 5E-1,\ 5E-2,\ 5E-3, 5E-4 and 5E-5. }
	\label{fig:fig2}
\end{figure}

%In this section, we test the RK-RFM for multiphase flow problem \eqref{eq:multiphasepde}  by  numerically simulating 
We first test the convergence of our method. The exact solution of the multiphase flow problem \eqref{eq:multiphasepde} is given by 
\begin{equation}\label{eq:exactsolution}
	\left\{
	\begin{aligned}
		&\phi_1(\boldsymbol{x},t)=\sin(x_1)\sin(x_2)\exp(-t),\\
		&\phi_2(\boldsymbol{x},t)=\cos(x_1)\cos(x_2)\exp(-t),\\
	\end{aligned}
	\right.
\end{equation}
for $(\boldsymbol{x},t)= (0,2\pi)^2\times(0,T]$ with initial condition
\begin{equation}\label{eq:exactsolutioninitial}
	\left\{
	\begin{aligned}
		&\phi_1(\boldsymbol{x},0)=\sin(x_1)\sin(x_2),\\
		&\phi_2(\boldsymbol{x},0)=\cos(x_1)\cos(x_2),\\
	\end{aligned}
	\right.
\end{equation}
% and periodic boundary condition
% \begin{equation}\label{eq:exactsolutionboundary}
% 	\left\{
% 	\begin{aligned}
% 		&\boldsymbol{\phi}(0,x_2,t)=\boldsymbol{\phi}(2\pi,x_2,t),\\
% 		&\boldsymbol{\phi}(x_1,0,t)=\boldsymbol{\phi}(x_1,2\pi,t),\\
% 	\end{aligned}
%         % \quad \boldsymbol{x}\in\partial\Omega
% 	\right.
% \end{equation}
and Dirichlet boundary condition on $\partial\Omega$
\begin{equation}\label{eq:exactsolutionboundary}
	\left\{
	\begin{aligned}
		&\phi_1(\boldsymbol{x},t)=\sin(x_1)\sin(x_2)\exp(-t),\\
		&\phi_2(\boldsymbol{x},t)=\cos(x_1)\cos(x_2)\exp(-t).
	\end{aligned}
%        \quad \boldsymbol{x}\in\partial\Omega
	\right.
\end{equation}
We set $T = 1.0$, $R_m=1.7$, $\gamma = \mu = \kappa =0.001$ and $R= \lambda =   \xi = \zeta = 1.0$.  
Algorithm \ref{alg:Algorithm 1} is employed to estimate $\boldsymbol{\phi}(\boldsymbol{x},t) = (\phi_1,\phi_2)^T$ with tanh and cos activation functions. 
The numerical results of the  errors for $\boldsymbol{\phi}$  with different  $\Delta t$ are shown in Table \ref{Table:1} and Figure \ref{fig:fig2}, where $Q_{test}$ is the number of test points. From Table \ref{Table:1} and Figure \ref{fig:fig2}, for the tanh activation function, it can be observed that when $\Delta t = 5E-1, 5E-2$ and $5E-3$, $e_1(\Delta t)< e_2(\Delta t)$,  and the relative $L^{\infty}$ and  $L^2$ errors achieve second-order convergence. When $\Delta t = 5E-4$, $e_1(\Delta t)\approx e_2(\Delta t)$,  and the relative $L^{\infty}$ and  $L^2$ errors reach the minimum in the five experiments. When $\Delta t = 5E-5$, $e_1(\Delta t)> e_2(\Delta t)$, the relative $L^{\infty}$ and  $L^2$ errors begin to increase. For the cos activation function, it can be observed that when $\Delta t = 5E-1, 5E-2, 5E-3$ and $5E-4$, $e_1(\Delta t)< e_2(\Delta t)$,  and the relative $L^{\infty}$ and  $L^2$ errors achieve second-order convergence. When $\Delta t = 5E-5$, $e_1(\Delta t)\approx e_2(\Delta t)$,  and the relative $L^{\infty}$ and  $L^2$ errors reach the minimum in the five experiments. 

\begin{figure}[htp]
	\centering 	 
		\includegraphics[scale=0.34]{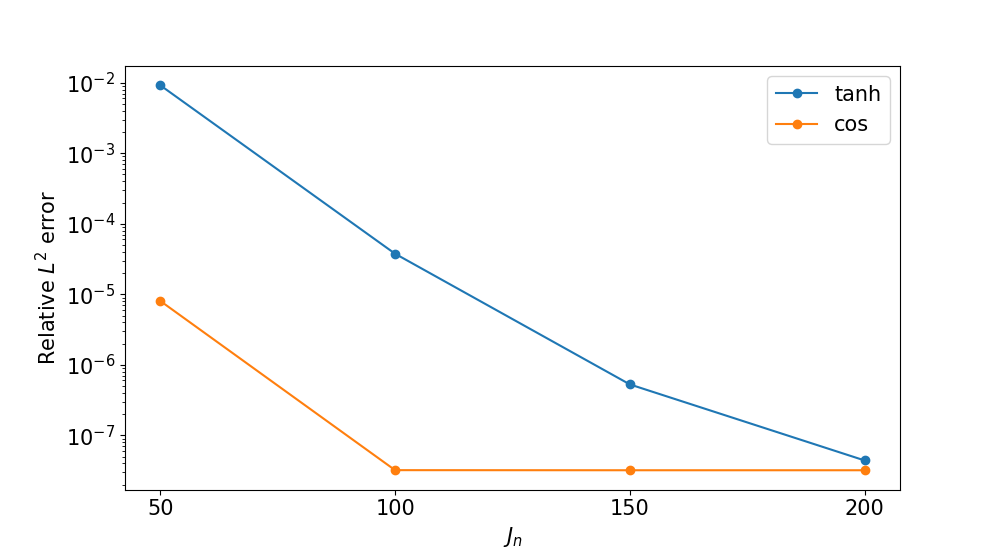}
	\caption{Relative $L^2$ error for \eqref{eq:exactsolution} obtained by the Algorithm \ref{alg:Algorithm 1} with $J_n=$ 50, 100,\ 150\ and 200. }
	\label{fig:fig3}
\end{figure}

\begin{figure}[htp]
	\centering 	 
		\includegraphics[scale=0.34]{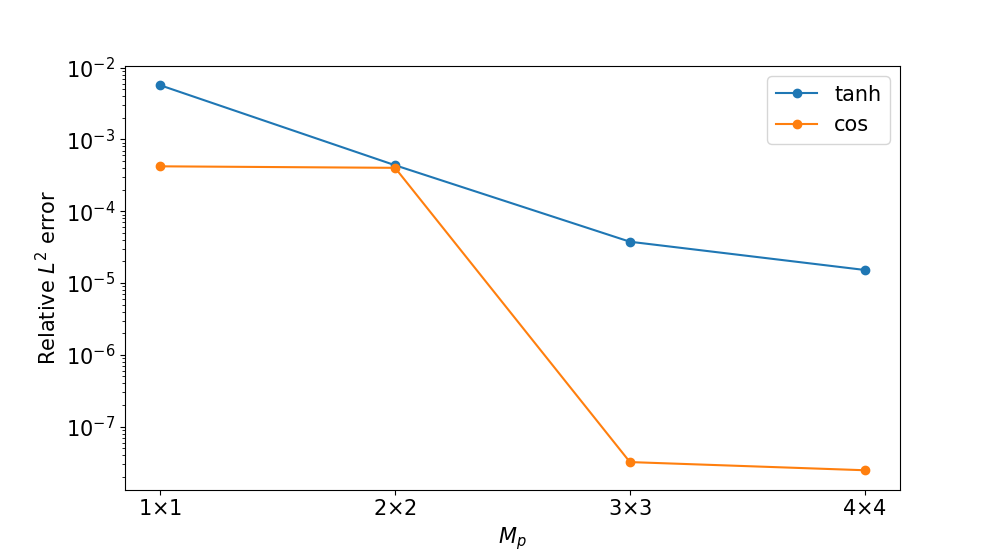}
	\caption{Relative $L^2$ error for \eqref{eq:exactsolution} obtained by the Algorithm \ref{alg:Algorithm 1} with $M_p=$ 1$\times$1,\ 2$\times$2,\ 3$\times$3 and 4$\times$4.}
	\label{fig:fig4}
\end{figure}

\begin{figure}[htp]
	\centering 	 
		\includegraphics[scale=0.34]{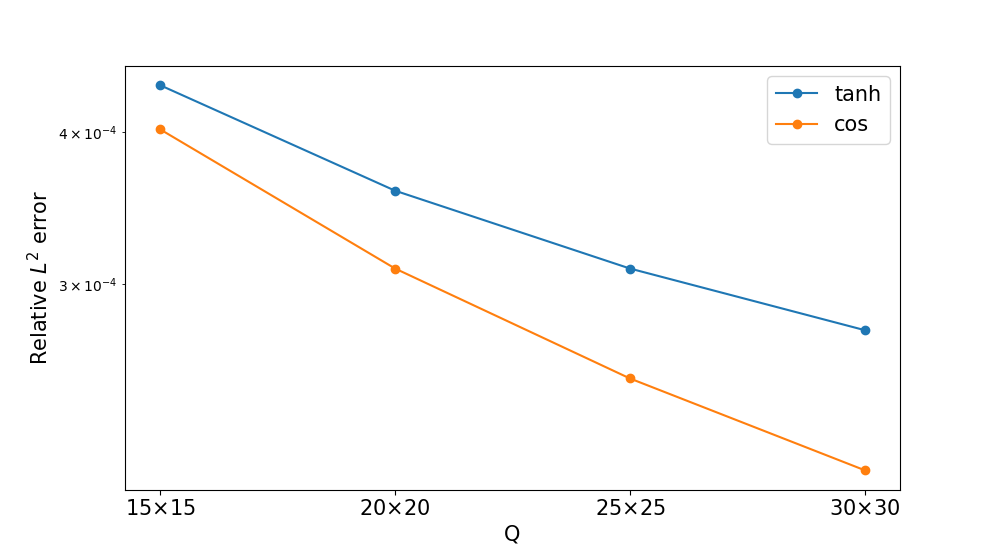}
	\caption{Relative $L^2$ error for \eqref{eq:exactsolution} obtained by the Algorithm \ref{alg:Algorithm 1} with $Q=$ 15$\times$15,\ 20$\times$20,\ 25$\times$25 and 30$\times$30. }
	\label{fig:fig5}
\end{figure}

In order to verify the universal approximation using our method, we respectively test the Algorithm \ref{alg:Algorithm 1} with different  $J_n$,  $M_P$ and $Q$ when $\Delta t = 5E-4$. The numerical results of the  errors for $\boldsymbol{\phi}$  are shown in Table \ref{Table:2}, Table \ref{Table:3}, Table \ref{Table:4}, Figure \ref{fig:fig3}, Figure \ref{fig:fig4} and Figure \ref{fig:fig5}, respectively. %From Table \ref{Table:2}, Table \ref{Table:3}, Table \ref{Table:4}, Figure \ref{fig:fig3}, Figure \ref{fig:fig4} and Figure \ref{fig:fig5}, 
From Table \ref{Table:2}, Table \ref{Table:3}, Figure \ref{fig:fig3} and Figure \ref{fig:fig4}, it can be observed that the relative $L^{\infty}$ and  $L^2$ errors decrease  as $J_n$ (number of the random feature functions on each subdomain) and $M_p$ (number of the subdomains) increase,  which is consistent with the Lemma \ref{lemma1}. From Table \ref{Table:4} and Figure \ref{fig:fig5}, it can be observed that the relative $L^{\infty}$ and  $L^2$ errors decrease  as $Q$ (number of collocation points in each subdomain)  increases, which  implies that increasing $Q$ can decrease $\delta$, which is defined in  \eqref{eq:definition1}.
% %and Figure \ref{fig:fig3}  depicts the training processes. 
% The relative $L^{2}$ errors and the training losses with different training steps are shown in Figure \ref{fig:fig3}. It is observed that these values decrease with parameter $\nu$ decreases. Similar phenomena will occur in the later experiments. Our method is not sensitive to the parameter $\nu$.
% The training time is 500s for each case, which is a acceptable cost. 
\begin{table*}[htp]
	\begin{center} 
 \caption{Relative $L^{\infty}$ and  $L^2$ errors for \eqref{eq:exactsolution} obtained by the Algorithm \ref{alg:Algorithm 1} with different $\Delta t$}
 %5E-1,\ 5E-2,\ 5E-3, 5E-4 and 5E-5. }
		\begin{tabular}{|l|l|l|l|l|l|l|l|l|l|l|}
			\hline
		$M_p$&$J_n$& $Q$&  $Q_{test}$& $\sigma$	& $\Delta t$ & 5E-1 & 5E-2  & 5E-3 & 5E-4 & 5E-5        \\
			\hline
	\multirow{6}{*}{3$\times$3}&\multirow{6}{*}{200}&\multirow{6}{*}{20$\times$20}&\multirow{6}{*}{40$\times$40} &\multirow{3}{*}{tanh}	&$\vert\vert e\vert\vert_{L^{\infty}} $  & 3.60E-2 & 3.53E-4 & 3.58E-6 & 1.45E-7 & 3.18E-7    \\
			\cline{6-11}
		& & & & &$\vert\vert e\vert\vert_{L^2} $  & 3.25E-2 & 3.15E-4 & 3.15E-6 & 3.70E-8  & 5.18E-8 \\
                \cline{6-11}
           & & &  & &time(s) & 0.75 & 5.28 & 50.65 & 500.75 & 5023.06   \\
           \cline{5-11}
        & & & &	\multirow{3}{*}{cos}&$\vert\vert e\vert\vert_{L^{\infty}} $  & 4.43E-2 & 3.58E-4 & 3.58E-6 & 3.79E-8 & 5.70E-10    \\
			\cline{6-11}
		& & & & &$\vert\vert e\vert\vert_{L^2} $  & 3.26E-2 & 3.15E-4 & 3.14E-6 & 3.18E-8  & 3.71E-10 \\
                \cline{6-11}
           & & &  & &time(s) & 0.81 & 4.98 & 47.46 & 470.18 & 4653.23   \\
			\hline
		\end{tabular}
		\label{Table:1}
	\end{center}
\end{table*}

% \begin{table}[htp]
% 	\begin{center} 
% 		\begin{tabular}{|l|l|l|l|l|l|l|l|l|l|}
% 			\hline
% 		$M_p$& $Q$&  $Q_{test}$	& $\Delta t$ & $J_n$ & 100 & 200  & 300 & 400       \\
% 			\hline
% 	\multirow{6}{*}{3$\times$3}&\multirow{6}{*}{26$\times$26}&\multirow{6}{*}{51$\times$51}&\multirow{3}{*}{5E-1}	&$\vert\vert e\vert\vert_{L^{\infty}} $  & 4.05E-2 & 4.04E-2 & 4.04E-2 & 4.04E-2    \\
% 			\cline{5-9}
% 		& & & &$\vert\vert e\vert\vert_{L^2} $  & 2.38E-2 & 2.38E-2 & 2.38E-2 & 2.38E-2  \\
%                 \cline{5-9}
%            & & &  &time(s) & 0.73 & 0.93 & 1.20 &  1.50  \\
% 			\cline{4-9}
% 	& & &\multirow{3}{*}{5E-5}	&$\vert\vert e\vert\vert_{L^{\infty}} $  & 1.13E-3 & 1.03E-6 & 1.79E-8 & 1.15E-8   \\
% 			\cline{5-9}
% 		& & & &$\vert\vert e\vert\vert_{L^2} $  & 4.97E-4 & 2.23E-7 & 6.27E-9 & 5.55E-9 \\
%                 \cline{5-9}
%            & & &  &time(s) & 4665.83 & 6470.20 & 8917.05 & 11648.99   \\
% 			\hline
% 		\end{tabular}
% 		\caption{Relative $L^{\infty}$ and  $L^2$ errors for \eqref{eq:exactsolution} performed by the Algorithm \ref{alg:Algorithm 1} with $\Delta t=$ 1E-1,\ 1E-2,\ 1E-3, 1E-4 and 1E-5. }\label{Table:2}
% 	\end{center}
% \end{table}

\begin{table*}[htp]
	\begin{center} 
      \caption{Relative $L^{\infty}$ and  $L^2$ errors for \eqref{eq:exactsolution} obtained by the Algorithm \ref{alg:Algorithm 1} with different $J_n$ }
		\begin{tabular}{|l|l|l|l|l|l|l|l|l|l|l|}
			\hline
		$M_p$& $Q$&  $Q_{test}$	& $\Delta t$ & $\sigma$	 & $J_n$ &50 & 100 & 150  & 200       \\
			\hline
	\multirow{6}{*}{3$\times$3}&\multirow{6}{*}{15$\times$15}&\multirow{6}{*}{30$\times$30}&\multirow{6}{*}{5E-4} &\multirow{3}{*}{tanh}	&$\vert\vert e\vert\vert_{L^{\infty}} $& 1.75E-2 & 1.01E-4 & 2.88E-6 & 3.38E-7  \\
			\cline{6-10}
		& & & & &$\vert\vert e\vert\vert_{L^2} $&9.21E-3  & 3.76E-5 & 5.27E-7 & 4.40E-8  \\
                \cline{6-10}
           & & & & &time(s)&224.44 & 271.97 & 301.70 & 370.18  \\
           \cline{5-10}
	& & & &\multirow{3}{*}{cos}	&$\vert\vert e\vert\vert_{L^{\infty}} $& 1.92E-5 & 4.02E-8 & 3.93E-8 & 3.91E-8  \\
			\cline{6-10}
		& & & & &$\vert\vert e\vert\vert_{L^2} $&8.06E-6  & 3.21E-8 & 3.20E-8 & 3.20E-8  \\
                \cline{6-10}
           & & & & &time(s)&202.81 & 278.36 & 319.03 & 344.35  \\
			\hline
		\end{tabular}
		\label{Table:2}
	\end{center}
\end{table*}

% \begin{table}[htp]
% 	\begin{center} 
% 		\begin{tabular}{|l|l|l|l|l|l|l|l|l|l|}
% 			\hline
% 		$J_n$& $Q$&  $Q_{test}$	& $\Delta t$ & $M_p$ & 1$\times$1 & 2$\times$2  & 3$\times$3 & 4$\times$4         \\
% 			\hline
% 	\multirow{6}{*}{200}&\multirow{6}{*}{26$\times$26}&\multirow{6}{*}{51$\times$51}&\multirow{3}{*}{5E-1}	&$\vert\vert e\vert\vert_{L^{\infty}} $  & 3.96E-2 & 3.88E-2 & 4.04E-2 & 3.85E-2    \\
% 			\cline{5-9}
% 		& & & &$\vert\vert e\vert\vert_{L^2} $  & 2.29E-2 & 2.28E-2 & 2.38E-2 & 2.27E-2  \\
%                 \cline{5-9}
%            & & &  &time(s) & 0.40 & 0.49 & 0.93 & 2.04   \\
% 			\cline{4-9}
% 	& & &\multirow{3}{*}{5E-5}	&$\vert\vert e\vert\vert_{L^{\infty}} $  & 4.60E-4 & 4.39E-5 & 1.03E-6 & 6.65E-7    \\
% 			\cline{5-9}
% 		& & & &$\vert\vert e\vert\vert_{L^2} $  & 1.44E-4 & 2.57E-5 & 2.23E-7 & 2.81E-7 \\
%                 \cline{5-9}
%            & & &  &time(s) & 1109.25 & 2685.99 & 6470.20 & 17579.80   \\
% 			\hline
% 		\end{tabular}
% 		\caption{Relative $L^{\infty}$ and  $L^2$ errors for \eqref{eq:exactsolution} performed by the Algorithm \ref{alg:Algorithm 1} with $\Delta t=$ 1E-1,\ 1E-2,\ 1E-3, 1E-4 and 1E-5. }\label{Table:3}
% 	\end{center}
% \end{table}

\begin{table*}[htp]
	\begin{center} 
     \caption{Relative $L^{\infty}$ and  $L^2$ errors for \eqref{eq:exactsolution} obtained by the Algorithm \ref{alg:Algorithm 1} with different $M_p$ }
		\begin{tabular}{|l|l|l|l|l|l|l|l|l|l|l|}
			\hline
		$J_n$& $Q$&  $Q_{test}$	& $\Delta t$ & $\sigma$ & $M_p$ & 1$\times$1 & 2$\times$2  & 3$\times$3 & 4$\times$4       \\
			\hline
	\multirow{6}{*}{100}&\multirow{6}{*}{15$\times$15}&\multirow{6}{*}{30$\times$30}&\multirow{6}{*}{5E-4} & \multirow{3}{*}{tanh}	&$\vert\vert e\vert\vert_{L^{\infty}} $  & 1.72E-2 & 9.00E-4 & 1.01E-4 & 4.84E-5   \\
			\cline{6-10}
		& & & & &$\vert\vert e\vert\vert_{L^2} $  & 5.68E-3 & 4.37E-4 & 3.76E-5 & 1.52E-5  \\
                \cline{6-10}
           & & & & &time(s) & 100.12 & 151.22 & 271.97 & 551.76   \\
           \cline{5-10}
	& & &  & \multirow{3}{*}{cos}	&$\vert\vert e\vert\vert_{L^{\infty}} $  & 7.61E-4 & 4.68E-4 & 4.02E-8 & 3.60E-8   \\
			\cline{6-10}
		& & & & &$\vert\vert e\vert\vert_{L^2} $  & 4.22E-4 & 4.02E-4 & 3.21E-8 & 2.47E-8 \\
                \cline{6-10}
           & & & & &time(s) & 95.07 & 140.55 & 278.36 & 513.32   \\
			\hline
		\end{tabular}
		\label{Table:3}
	\end{center}
\end{table*}

\begin{table*}[htp]
	\begin{center} 
     \caption{Relative $L^{\infty}$ and  $L^2$ errors for \eqref{eq:exactsolution} obtained by the Algorithm \ref{alg:Algorithm 1} with different $Q$ }
		\begin{tabular}{|l|l|l|l|l|l|l|l|l|l|l|}
			\hline
		\multirow{2}{*}{$M_p$}& \multirow{2}{*}{$J_n$}&  \multirow{2}{*}{$\Delta t$}  & \multirow{2}{*}{$\sigma$}& $Q$ & 15$\times$15 & 20$\times$20  & 25$\times$25 & 30$\times$30       \\
			\cline{5-9}
        & &  & & $Q_{test}$ & 30$\times$30 & 40$\times$40  & 50$\times$50 &  60$\times$60      \\
			\hline
	\multirow{6}{*}{2$\times$2}&\multirow{6}{*}{100}&\multirow{6}{*}{5E-4} & \multirow{3}{*}{tanh}	&$\vert\vert e\vert\vert_{L^{\infty}} $  & 9.00E-4 & 8.36E-4 & 7.81E-4 & 7.41E-4   \\
			\cline{5-9}
		& & & &$\vert\vert e\vert\vert_{L^2} $  & 4.37E-4 & 3.58E-4 & 3.09E-4 & 2.75E-4 \\
                \cline{5-9}
           & & &  &time(s) & 151.22 & 173.75 & 198.86 & 259.54  \\
           \cline{4-9}
	 & & & \multirow{3}{*}{cos}	&$\vert\vert e\vert\vert_{L^{\infty}} $  & 4.68E-4 & 3.68E-4 & 3.01E-4 & 2.55E-4   \\
			\cline{5-9}
		&  &  & &$\vert\vert e\vert\vert_{L^2} $  & 4.02E-4 & 3.09E-4 & 2.51E-4 & 2.11E-4 \\
                \cline{5-9}
           & &  & &time(s) & 140.55 & 161.65 & 188.70 & 250.59  \\
			\hline
		\end{tabular}
		\label{Table:4}
	\end{center}
\end{table*}
\subsection{Multiphase Flow Problem of Cells}
\label{sec06:sec02}
\begin{figure*}[htp]
	\centering 
	\subfigure[t=0]{		 
		\includegraphics[scale=0.31]{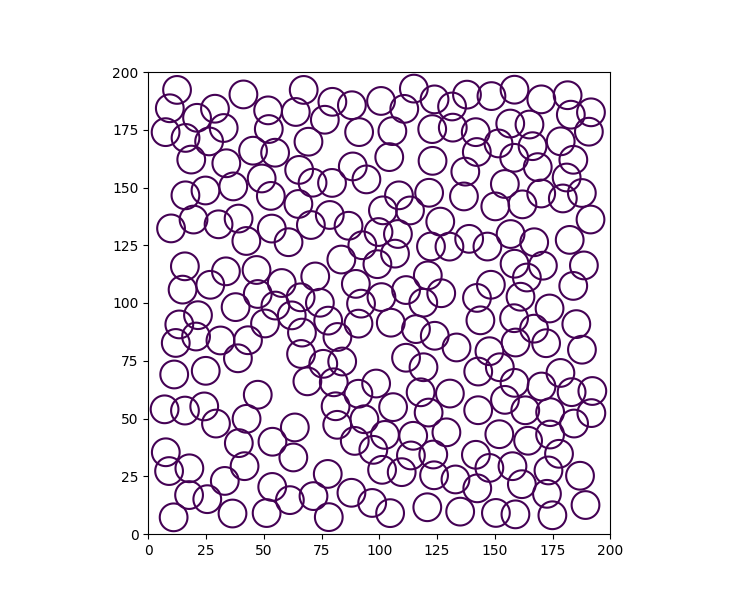}}
            \subfigure[t=10]{		 
            \includegraphics[scale=0.31]{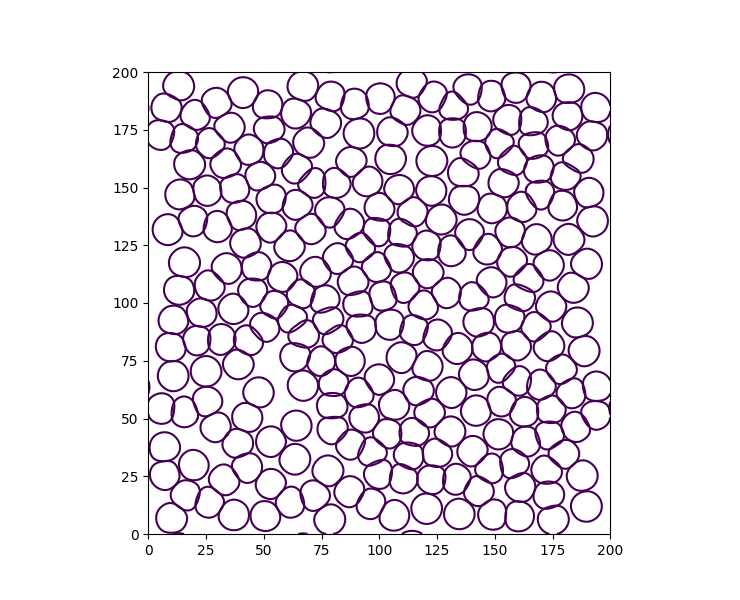}}
            \subfigure[t=30]{		 
            \includegraphics[scale=0.31]{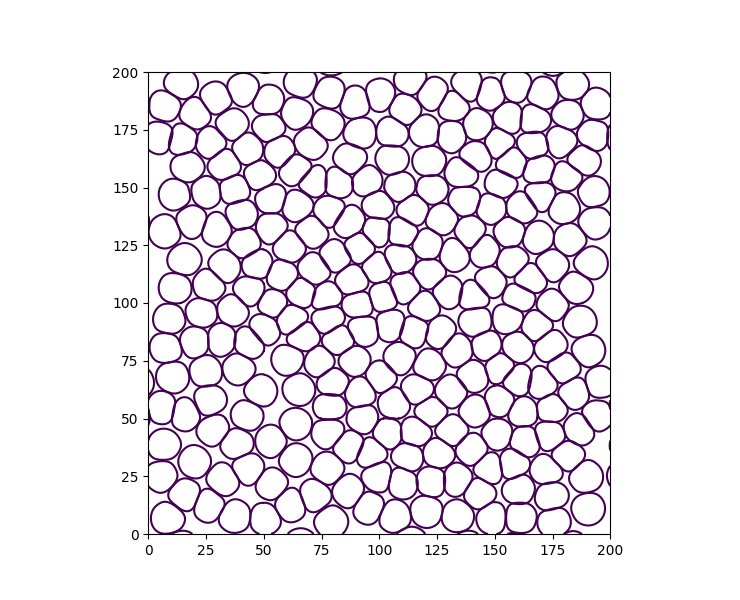}}
            \subfigure[t=50]{		 
            \includegraphics[scale=0.31]{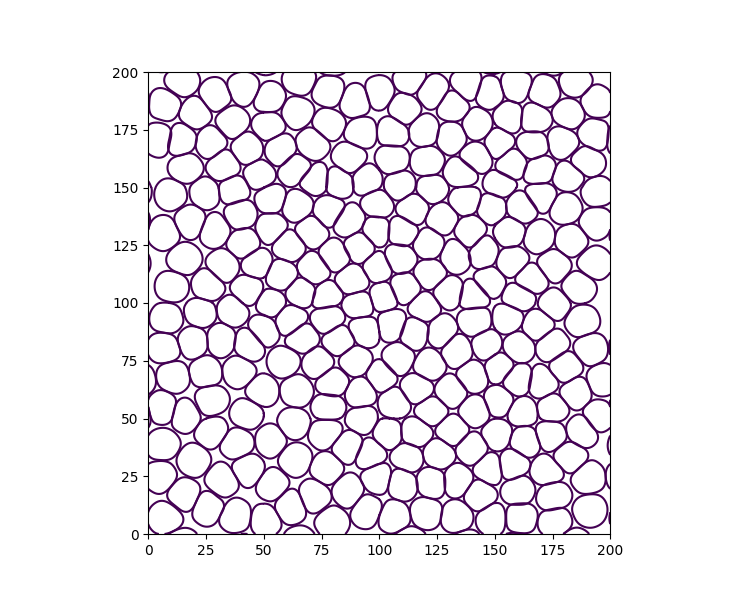}}
            \subfigure[t=70]{		 
            \includegraphics[scale=0.31]{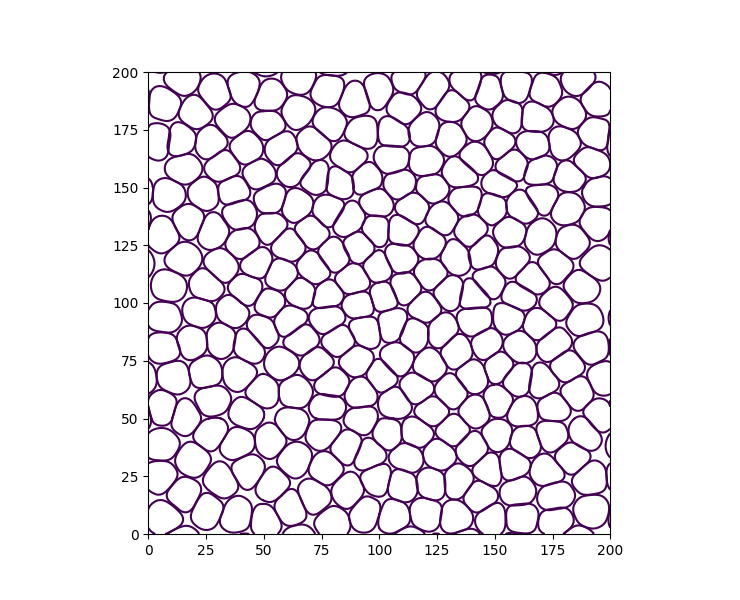}}
            \subfigure[t=100]{		 
            \includegraphics[scale=0.31]{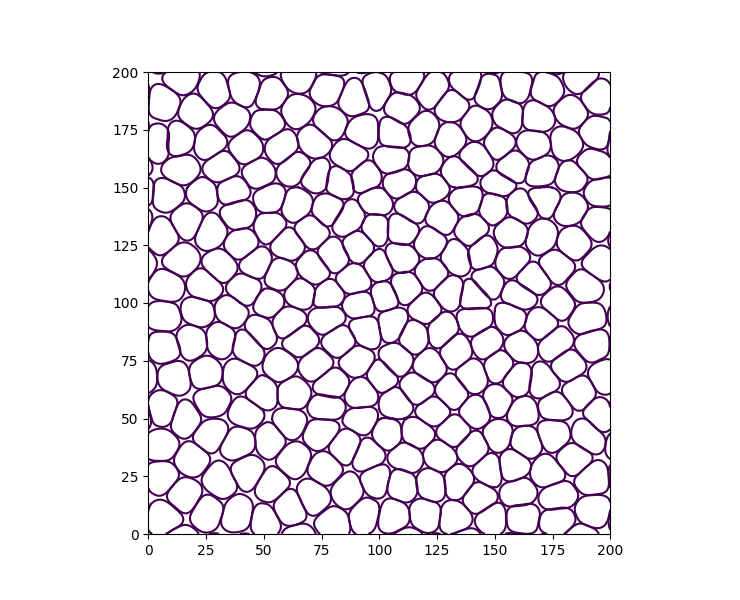}}
	\caption{The evolution of the cells at $t = 0, 10, 30, 50, 70, 100$}
	\label{fig:fig6}
\end{figure*}

Now, we use the RK-RFM to solve the  multiphase flow problem of cells \eqref{eq:multiphasepde} for $(\boldsymbol{x},t)= (0,200)^2\times(0,T]$ with 
initial conditions
\begin{equation}\label{eq:initialcell}
	\begin{aligned}
		&\phi_{i}(\boldsymbol{x}, 0)=\frac{1+\tanh \left(\frac{r-\vert\vert \boldsymbol{x}-C_i\vert\vert_{L^2}}{\lambda}\right)}{2},i=1, \ldots, d_{\phi},
	\end{aligned}
        % \quad \boldsymbol{x}\in\partial\Omega
\end{equation} 
and periodic boundary condition
\begin{equation}\label{eq:periodicboundary}
	\left\{
	\begin{aligned}
		&\boldsymbol{\phi}(0,x_2,t)=\boldsymbol{\phi}(200,x_2,t),\\
		&\boldsymbol{\phi}(x_1,0,t)=\boldsymbol{\phi}(x_1,200,t),\\
	\end{aligned}
        % \quad \boldsymbol{x}\in\partial\Omega
	\right.
\end{equation} 
where $r$ is the initial radiu and  $C_i$ is the initial center of each cell.
Algorithm \ref{alg:Algorithm 1} is employed to simulate the dynamic of $\boldsymbol{\phi}(\boldsymbol{x},t)$. We set $T = 100.0$, $R_m=5.0$, $d_{\phi} = 240$, $r=6.0$, $\gamma = 0.01$, $\mu=3$, $\kappa = 0.1$, $\lambda = 2.5$, $R=8.0$, $\mu=3$, $\xi=2$ and $\zeta = 0.005$. The parameters of the Algorithm \ref{alg:Algorithm 1} are chosen to be $\Delta t=0.1$,  $M_p = 4\times 4$, $Q=40\times40$, $J_n=600$ and $Q_{test}=200\times200$. The $\tanh$ function is taken as the activation function. Initially, we randomly generate cells with radius $r$ and center $C_i$ in $\Omega$. The area of the cells $\pi r^2$ is smaller than the target area of $\pi R^2$. We allow them to partially overlap and relax during a period of inactivity.  This simulation costs 1h 37min. The dynamical behavior of the cells is visually shown in 
Figure \ref{fig:fig6}, and the configuration is in a steady state after a long time evolution.

\begin{figure*}[htp]
	\centering
    \subfigure[]{
		\includegraphics[scale=0.34]{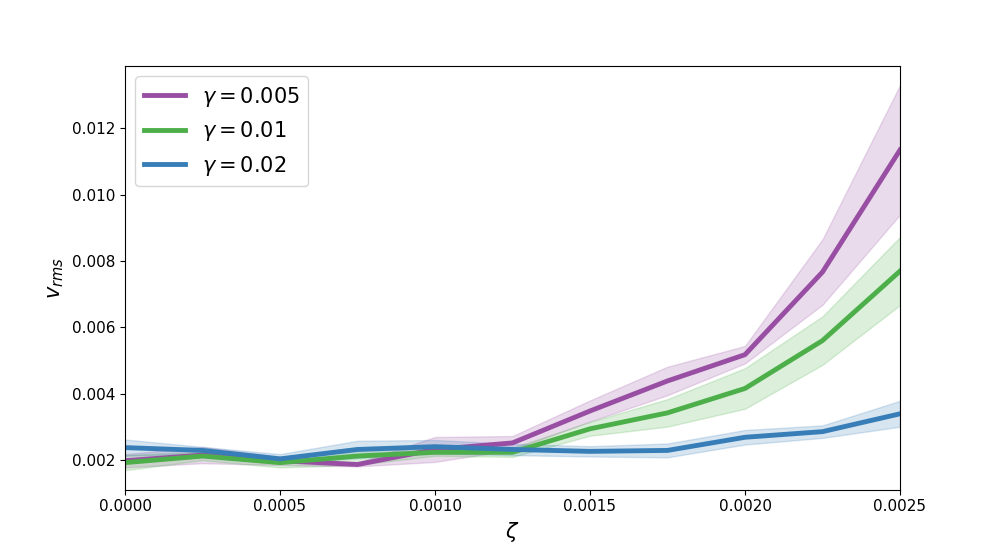}}
	\subfigure[]{		 
      \includegraphics[scale=0.34]{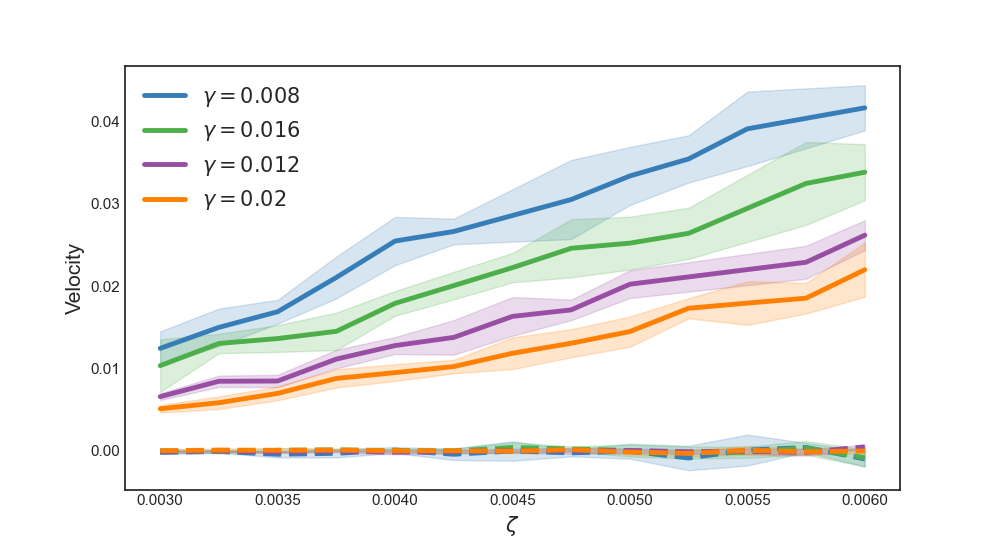}}
      \subfigure[]{
		\includegraphics[scale=0.34]{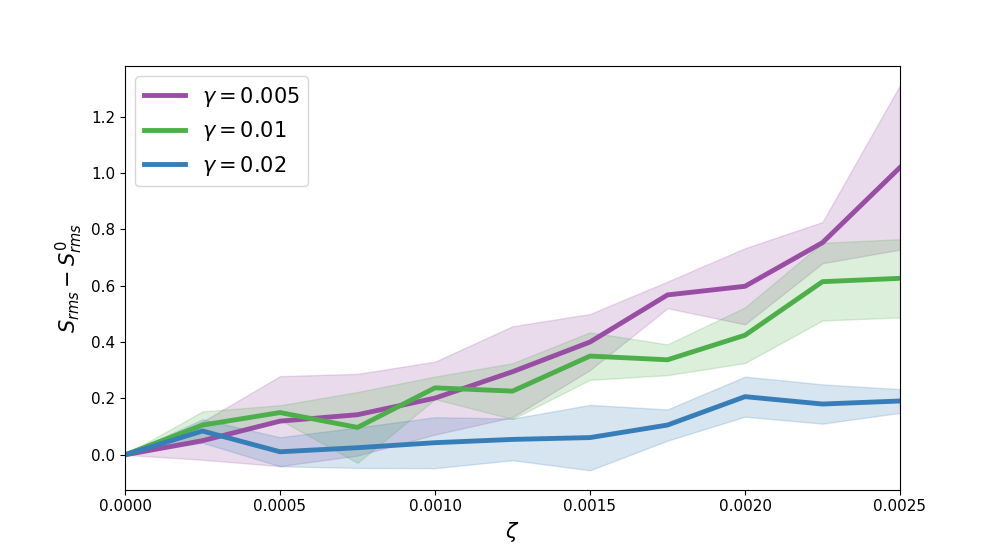}}
  \subfigure[]{		 
		\includegraphics[scale=0.34]{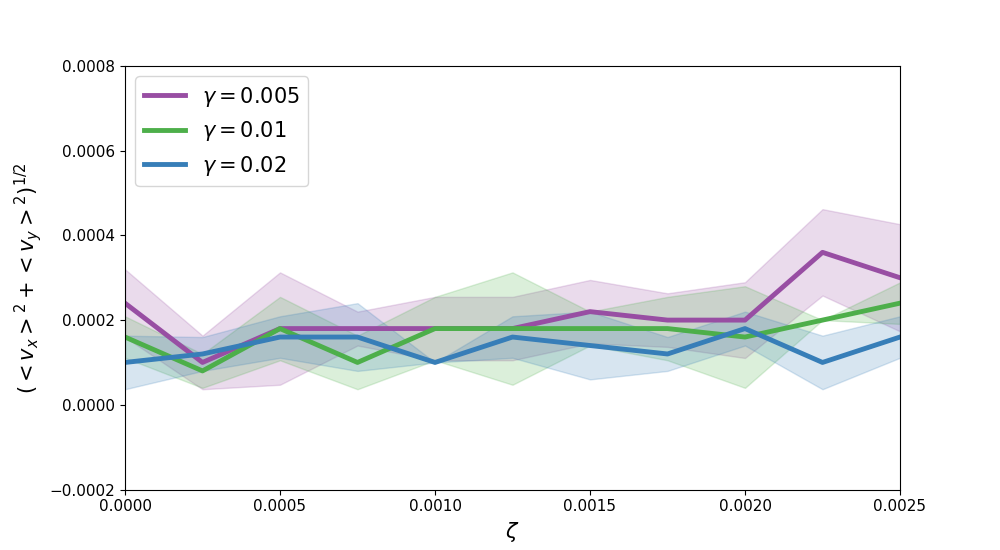}}
  \subfigure[]{		 
		\includegraphics[scale=0.34]{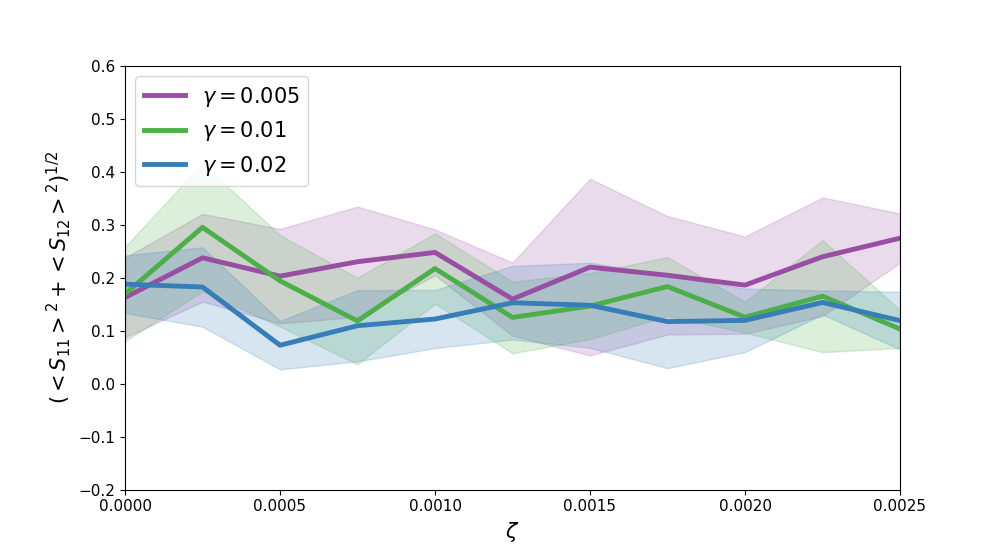}}
      \subfigure[]{		 
		\includegraphics[scale=0.34]{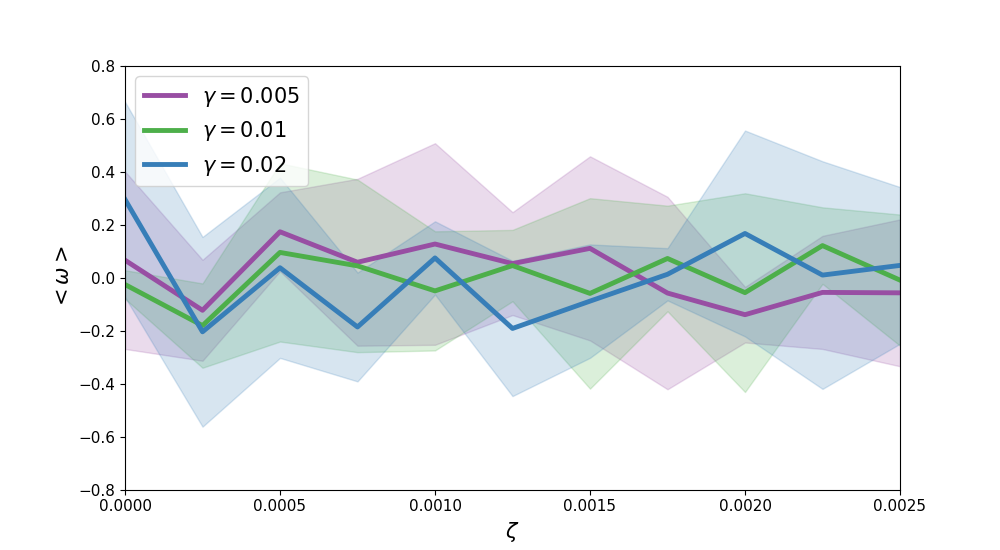}}
%    \subfigure[]{		 
		\caption{The values of of $v_{\mathrm{rms}}$, $\left\langle v\right\rangle$, $S_{\mathrm{rms}}-S_{\mathrm{rms}}^0$, $\left(\left\langle v_x\right\rangle^2+\left\langle v_y\right\rangle^2\right)^{1 / 2}$, $\left(\left\langle S_{11}\right\rangle^2+\left\langle S_{12}\right\rangle^2\right)^{1 / 2}$ and $\left\langle\omega\right\rangle$ with different $\gamma$ and $\zeta$  at $t = 1000.0$. (b) $v_{\mathrm{rms}}$ (plain lines) and $\left\langle v\right\rangle$ (dashed lines). Mean $\pm$ std from five simulations.}
	\label{fig:fig7}
\end{figure*}

In order to further verify  the cells evolution obtained by Algorithm \ref{alg:Algorithm 1}, we define the tissue
velocity as $v=\sum_{i=1}^{d_{\phi}} \phi_i v_i$.
Then we measure the root-mean-square velocity $v_{\mathrm{rms}}=\left(\left\langle v^2\right\rangle\right)^{1/2}$ and  nematic order $S_{\mathrm{rms}}=\left(\left\langle\operatorname{det} S^2\right\rangle\right)^{1/2}=\left(\left\langle S_{11}^2+S_{12}^2\right\rangle\right)^{1/2}$,  the square roots of the sums of the squares of the means of the individual components of the velocity $\left(\left\langle v_x\right\rangle^2+\left\langle v_y\right\rangle^2\right)^{1 / 2}$ and order $\left(\left\langle S_{11}\right\rangle^2+\left\langle S_{12}\right\rangle^2\right)^{1 / 2}$, the mean
angle of the director $\left\langle\omega\right\rangle=\left\langle\arctan \left(S_{12} / S_{11}\right)\right\rangle$ and the mean velocity $\left\langle v\right\rangle$ separately with different $\gamma$ and $\zeta$. Specifically, we set $\Omega = (0,100)^2$, $T = 1000.0$, $R_m=5.0$, $d_{\phi} = 60$, $r=6.0$,  $\mu=3$, $\kappa = 0.1$, $\lambda = 2.5$, $R=8.0$, $\mu=3$ and $\xi=2$. The parameters of the Algorithm \ref{alg:Algorithm 1} are chosen to be $\Delta t=1.0$,  $M_p = 3\times 3$, $Q=100\times100$, $J_n=600$, $Q_{test}=200\times200$ with the initial condition \eqref{eq:initialcell} and periodic boundary conditions. The $\tanh$ function is taken as the activation function. The numerical results  are shown in Figure \ref{fig:fig7}. %The changing trends of the values of $v_{\mathrm{rms}}$, $S_{\mathrm{rms}}$, $\left(\left\langle v_x\right\rangle^2+\left\langle v_y\right\rangle^2\right)^{1 / 2}$, $\left(\left\langle S_{11}\right\rangle^2+\left\langle S_{12}\right\rangle^2\right)^{1 / 2}$, $\left\langle\omega\right\rangle$ and $\left\langle v\right\rangle$ are consistent with the results in \cite{mueller2019emergence}, which implies that the evolution of cells can be reasonably simulated by our method.
It is observed that the values $v_{\mathrm{rms}}$ and $S_{\mathrm{rms}}-S^0_{\mathrm{rms}}$ develop nonzero values with increasing activity $\zeta$, where $S^0_{\mathrm{rms}} = S_{\mathrm{rms}}|_{\zeta =0}$. The values $\left(\left\langle v_x\right\rangle^2+\left\langle v_y\right\rangle^2\right)^{1 / 2}$ and $\left(\left\langle S_{11}\right\rangle^2+\left\langle S_{12}\right\rangle^2\right)^{1 / 2}$ do not exhibit a significant rise with increasing $\zeta$. The mean angle $\left\langle w \right \rangle$ is uniformly distributed in $[-\pi,\pi]$ so that its average value is almost zero. The changing trends of the values shown in Figure \ref{fig:fig7} are consistent with the results in  reference \cite{mueller2019emergence}.

\begin{figure*}[htp]
	\centering
    \subfigure[]{
		\includegraphics[scale=0.34]{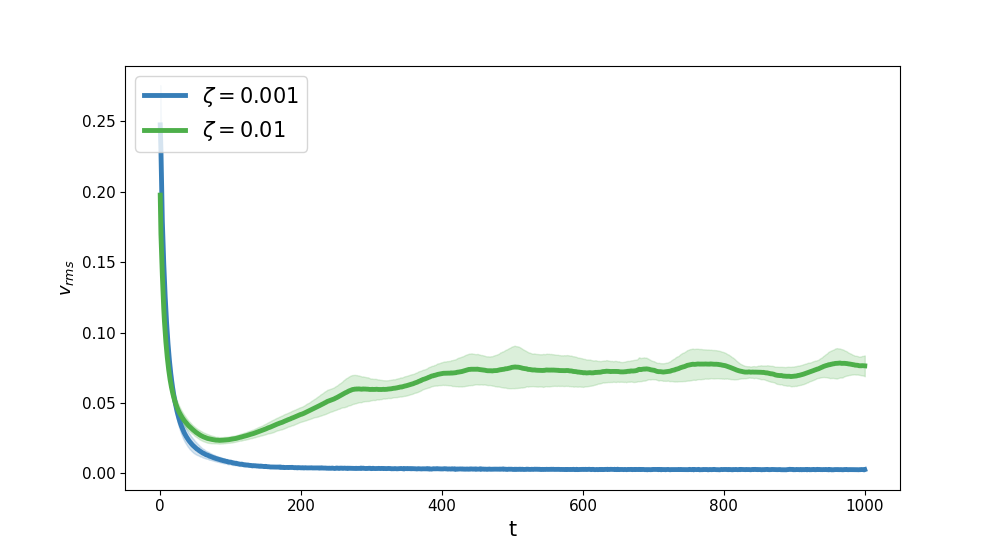}}
	\subfigure[]{		 
		\includegraphics[scale=0.34]{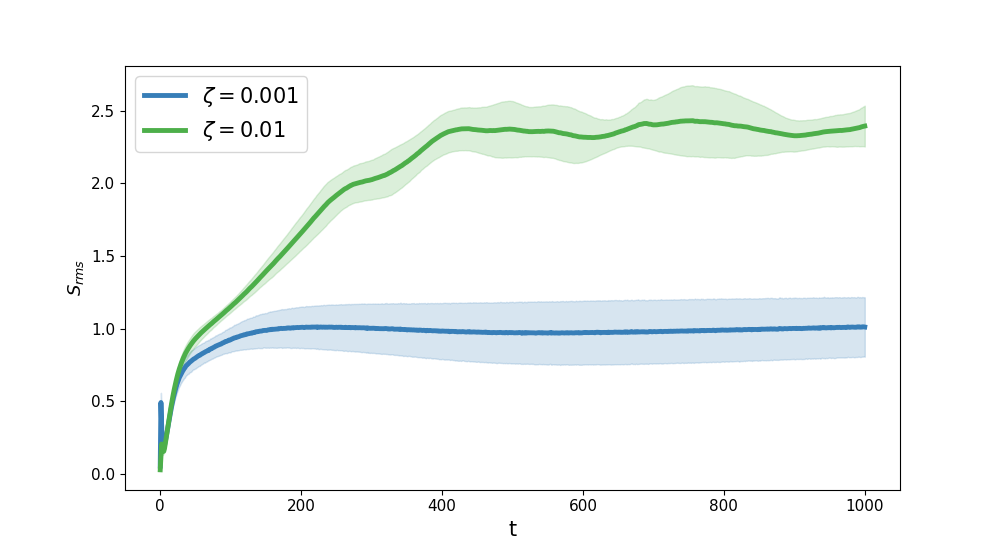}}
	\caption{The temporal evolution of  $v_{\mathrm{rms}}$ and  $S_{\mathrm{rms}}$ below (blue) and above (green) the activity threshold. Mean $\pm$ std from five simulations.}
	\label{fig:fig8}
\end{figure*}

\begin{figure*}[htp]
	\centering
    \subfigure[]{
		\includegraphics[scale=0.34]{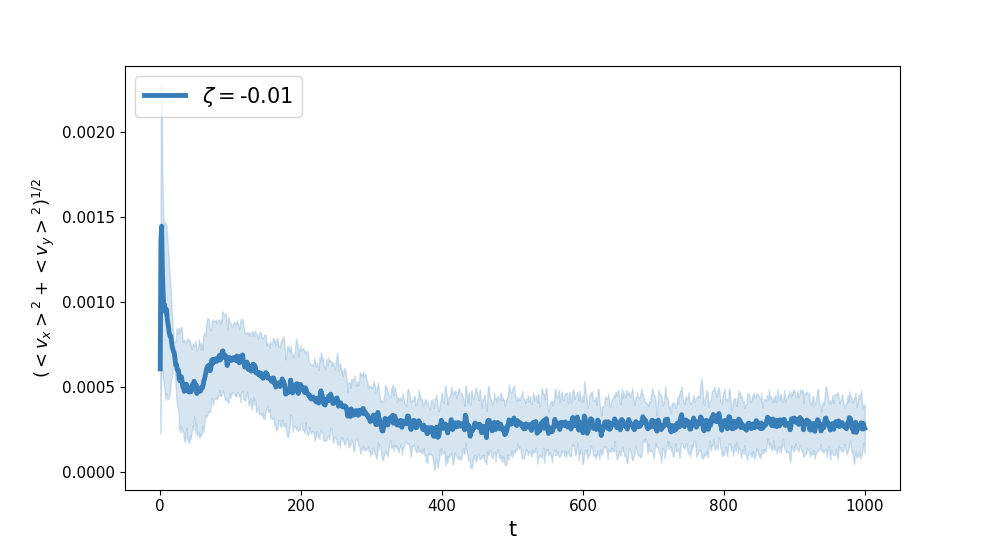}}
	\subfigure[]{		 
		\includegraphics[scale=0.34]{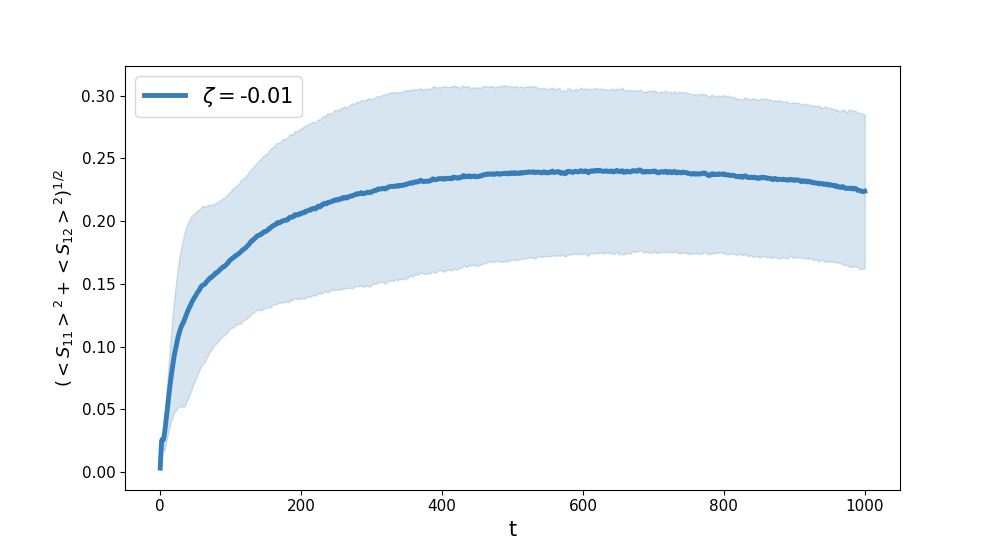}}
	\caption{The temporal evolution of  $\left(\left\langle v_x\right\rangle^2+\left\langle v_y\right\rangle^2\right)^{1 / 2}$ and  $\left(\left\langle S_{11}\right\rangle^2+\left\langle S_{12}\right\rangle^2\right)^{1 / 2}$ with  $\zeta = -0.01$. Mean $\pm$ std from five simulations.}
	\label{fig:fig9}
\end{figure*}

Finally, we investigate the temporal evolution of  $v_{\mathrm{rms}}$ and  $S_{\mathrm{rms}}$ below or above the activity threshold. 
The activity threshold  will emerge because cells must exert a strong enough push (or pull) to induce a sufficient deformation in their neighboring cells \cite{mueller2019emergence}. In addition, we investigate the temporal evolution of  $\left(\left\langle v_x\right\rangle^2+\left\langle v_y\right\rangle^2\right)^{1 / 2}$ and  $\left(\left\langle S_{11}\right\rangle^2+\left\langle S_{12}\right\rangle^2\right)^{1 / 2}$ with  $\zeta< 0$. Parameters and settings are the same as the example above. The temporal changes in $v_{\mathrm{rms}}$ and  $S_{\mathrm{rms}}$ below and above the activity threshold are shown are Figure \ref{fig:fig8}. In both cases, the initial relaxation from the base setup is noticeable. %The numerical results  are shown in Figure \ref{fig:fig8} and Figure \ref{fig:fig9}. 
%From the Figure \ref{fig:fig8}, 
 It is observed that the system reverts to equilibrium when the activity level is below the threshold. Conversely, when the activity level is above the threshold, the system can maintain large flows and exhibits a continuous increase in its order. From Figure \ref{fig:fig9}, it can be seen that the system is stable when $\zeta< 0$ which implies that the interactions at the cell interfaces tend to restore the stable state. The fluctuating patterns of the data depicted in Figure \ref{fig:fig9} align with the results in reference \cite{mueller2019emergence}. Since the system is initialized by creating many round cells defined in \eqref{eq:initialcell}, which is different from the initial setting in \cite{mueller2019emergence}, the value of nematic order increases for some initial time and attains
its steady state.

\section{Conclusions and remarks}
\label{sec05}
In this work, we present the Runge-Kutta random feature method (RK-RFM), a novel method designed to address the strongly nonlinear coupled multiphase flow problems of cells, while maintaining high numerical accuracy with reduced computational resources. The RK-RFM  employs the  RFM in space and the explicit RK method in time, achieving high accuracy in both space and time. This algorithm is mesh-free, making it suitable for complex geometric configurations, and is effective in solving PDEs with strong nonlinearity. Furthermore, we choose to use the same random feature functions for different cells and manually derive the necessary neural network derivative functions before the calculation begins, instead of relying on automatic differentiation, which significantly accelerates computation and conserves computational resources. We also provide error estimates for the RK-RFM and present some numerical experiments to validate our method and the application of RK-RFM in the multiphase flow problem of cells. In all numerical results, our method is verified to be both stable and accurate. In the future work, we will explore the potential of the RK-RFM in other types of PDEs and broader applications.

\begin{acknowledgments}
This research is partially supported by the National key R \& D Program of China (No.2022YFE03040002) and the National Natural Science Foundation of China ( No.12371434).
\end{acknowledgments}

\section*{Data Availability Statement}
The data that support the findings of this study are available from the corresponding author upon reasonable request.

\appendix

\section{}
\label{app:A}
To better illustrate how to obtain the optimal $\boldsymbol{U}_{t_{k+1}}$ by the linear least-squares method, we consider cases with $d=2$ and $d_{\phi}=1$. The arrangement of subdomains and the distribution of collocation points in subdomain $\Omega_{n}$ are set according to Figure \ref{fig:fig13}.
Now, we have
\begin{equation*}
	\begin{aligned}
\boldsymbol{x}^n=\{x^n_{1,1},\cdots,x^n_{1,Q_x},\cdots,x^n_{Q_y,1},\cdots,x^n_{Q_y,Q_x}\}.
\end{aligned}
\end{equation*}
% where Qx and Qy are the number of points in the x
% and y directions on the n-th sub-domain.
The coefficient matrix corresponding to $\tilde{\boldsymbol{\phi}}_{t_{k+1}}(\boldsymbol{x})$ in the subdomain $\Omega_{n}$ is defined as follows
\begin{equation*}
    \begin{aligned}
\boldsymbol{A}_n(\boldsymbol{x}^n)=\left[\begin{array}{ccc}\varphi_{n 1, t_{k+1}}(x^n_{1,1}) & \cdots & \varphi_{n J_n, t_{k+1}}(x^n_{1,1}) \\ \vdots & \cdots & \vdots \\\varphi_{n 1, t_{k+1}}(x^n_{Q_y,Q_x}) & \cdots & \varphi_{n J_n, t_{k+1}}(x^n_{Q_y,Q_x})\end{array}\right]_{Q \times J_n}.
    \end{aligned}
\end{equation*}
Then, we  get the coefficient matrix corresponding to $\tilde{\boldsymbol{\phi}}_{t_{k+1}}(\boldsymbol{x})$  in the entire domain $\Omega$
\begin{equation*}
    \begin{aligned}
\boldsymbol{A}=diag(\boldsymbol{A}_1(\boldsymbol{x}^1),\cdots,\boldsymbol{A}_{M_p}(\boldsymbol{x}^{M_p}))_{M_pQ \times M_pJ_n}.
\end{aligned}
\end{equation*}

Let \begin{equation*}
    \begin{aligned}
    &\boldsymbol{x}^n_{,1}=\{\boldsymbol{x}^n_{1,1},\cdots,\boldsymbol{x}^n_{Q_y,1}\},\
    \boldsymbol{x}^n_{,Q_x}=\{\boldsymbol{x}^n_{1,Q_x},\cdots,\boldsymbol{x}^n_{Q_y,Q_x}\},\\ 
    &\boldsymbol{x}^n_{1,}=\{\boldsymbol{x}^n_{1,1},\cdots,\boldsymbol{x}^n_{1,Q_x}\},\ \boldsymbol{x}^n_{Q_y,}=\{\boldsymbol{x}^n_{Q_y,1},\cdots,\boldsymbol{x}^n_{Q_y,Q_x}\},
    \end{aligned}
\end{equation*}
we have the coefficient matrices of the boundary conditions, where the boundary is perpendicular to the $x$-axis or $y$-axis in the subdomain $\Omega_{n}$, 
\begin{equation*}
    \begin{aligned}
&\boldsymbol{B}_{x,n,l}=[\boldsymbol{0}_{Q_y \times (n-1)J_n}, \mathcal{B}\boldsymbol{A}_n(\boldsymbol{x}^n_{,1}),\boldsymbol{0}_{Q_y \times (M_p-n)J_n}]_{Q_y \times M_pJ_n}, \\
 &\boldsymbol{B}_{x,n,r}=[\boldsymbol{0}_{Q_y \times (n-1)J_n}, \mathcal{B}\boldsymbol{A}_n(\boldsymbol{x}^n_{,Q_x}),\boldsymbol{0}_{Q_y \times (M_p-n)J_n}]_{Q_y \times M_pJ_n},\\
&\boldsymbol{B}_{y,n,d}=[\boldsymbol{0}_{Q_x \times (n-1)J_n}, \mathcal{B}\boldsymbol{A}_{n}(\boldsymbol{x}^{n}_{1,}),\boldsymbol{0}_{Q_x \times (M_p-n)J_n}]_{Q_x \times M_pJ_n}, \\
&\boldsymbol{B}_{y,n,u}=[\boldsymbol{0}_{Q_x \times (n-1)J_n}, \mathcal{B}\boldsymbol{A}_{n}(\boldsymbol{x}^{n}_{Q_y,}),\boldsymbol{0}_{Q_x \times (M_p-n)J_n}]_{Q_x \times M_pJ_n}.
\end{aligned}
\end{equation*}
where the subscripts $l$, $r$, $d$, and $u$ denote the left boundary, right boundary, lower boundary, and upper boundary, respectively.
Therefore, we obtain
\begin{equation*}
    \begin{aligned}
\boldsymbol{B}_x=&\left[\boldsymbol{B}_{x,1,l}^T,\cdots,\boldsymbol{B}_{x,N_y,l}^T,\right.\\
&\left.\boldsymbol{B}_{x,(N_x-1)N_y+1,r}^T,\cdots,\boldsymbol{B}_{x,N_xN_y,r}^T\right]^T_{2N_yQ_y \times M_pJ_n},\\
\boldsymbol{B}_y=&\left[\boldsymbol{B}_{y,1,d}^T,\boldsymbol{B}_{y,N_y,u}^T,\cdots,\right.\\
&\left.\boldsymbol{B}_{y,(N_x-1)N_y+1,d}^T,\boldsymbol{B}_{y,N_xN_y,u}^T\right]^T_{2N_xQ_x \times M_pJ_n},
\end{aligned}
\end{equation*}
in the whole domain $\Omega$.

In the subdomain  $\Omega_{n}$,  the coefficient matrix of the $C^0$ continuity condition is as follows, where the  interface is perpendicular to the $x$-axis or $y$-axis,  
\begin{equation*}
    \begin{aligned}
\boldsymbol{C}_{x,n}^0=&[\boldsymbol{0}_{Q_y \times (n-1)J_n}, \boldsymbol{A}_n(\boldsymbol{x}^n_{,Q_x}),\boldsymbol{0}_{Q_y \times N_yJ_n},\\
&-\boldsymbol{A}_{n+Ny}(\boldsymbol{x}^{n+Ny}_{,1}),\boldsymbol{0}_{Q_y \times (M_p-n-N_y)J_n}]_{Q_y \times M_pJ_n},\\
\boldsymbol{C}_{y,n}^0=&[\boldsymbol{0}_{Q_x \times (n-1)J_n}, \boldsymbol{A}_n(\boldsymbol{x}^n_{Q_y,}),\\
&-\boldsymbol{A}_{n+1}(\boldsymbol{x}^{n+1}_{1,}),\boldsymbol{0}_{Q_x \times (M_p-n-1)J_n}]_{Q_x \times M_pJ_n},
\end{aligned}
\end{equation*}
then we obtain that
\begin{equation*}
    \begin{aligned}
\boldsymbol{C}_{x}^0=& [(\boldsymbol{C}_{x,1}^0)^T,\cdots,(\boldsymbol{C}_{x,N_y(N_x-1)}^0)^T]^T_{N_y(N_x-1)Q_y \times M_pJ_n},\\
\boldsymbol{C}_{y}^0=& [(\boldsymbol{C}_{y,1}^0)^T,\cdots,(\boldsymbol{C}_{y,N_y-1}^0)^T,\cdots,\\
&(\boldsymbol{C}_{y,N_y(N_x-1)+1}^0)^T,\cdots,(\boldsymbol{C}_{y,N_yN_x-1}^0)^T]^T_{N_x(N_y-1)Q_x \times M_pJ_n}.
\end{aligned}
\end{equation*}

Similarly, we have the following coefficient matrices of the $C^{1}$ continuity condition
\begin{equation*}
    \begin{aligned}
&\boldsymbol{C}_{x}^0,\  \boldsymbol{C}_{x}^{1}=\nabla_x\boldsymbol{C}_{x}^0,\\ &\boldsymbol{C}_{y}^0,\  \boldsymbol{C}_{y}^{1}=\nabla_y\boldsymbol{C}_{y}^0.
\end{aligned}
\end{equation*}

The right-hand side term corresponding to $\boldsymbol{A}_n(\boldsymbol{x}^n)$ is defined as
\begin{equation*}
    \begin{aligned}
\boldsymbol{f}_{A_n}(\boldsymbol{x}^n) = [\boldsymbol{\phi}_{t_{k+1}}(x^n_{1,1}),\cdots,\boldsymbol{\phi}_{t_{k+1}}(x^n_{Q_y,Q_x})]_{1\times Q},
\end{aligned}
\end{equation*}
therefore we have
\begin{equation*}
    \begin{aligned}
    \boldsymbol{f}_A = [\boldsymbol{f}_{A_1}(\boldsymbol{x}^1),\cdots,\boldsymbol{f}_{A_{M_p}}(\boldsymbol{x}^{M_p})]^T_{M_pQ \times 1},
    \end{aligned}
\end{equation*}
related to $A$.

Also, we have
\begin{equation*}
    \begin{aligned}
\boldsymbol{f}_{B_x} =& [\boldsymbol{g}(\boldsymbol{x}^1_{,1},t_{k+1}),\cdots,\boldsymbol{g}(\boldsymbol{x}^{N_y}_{,1},t_{k+1}),\\
&\boldsymbol{g}(\boldsymbol{x}^{(N_x-1)N_y+1}_{,Q_x},t_{k+1}),\cdots,\boldsymbol{g}(\boldsymbol{x}^{N_xN_y}_{,Q_x},t_{k+1})]^T_{N_yQ_y \times 1},\\
\boldsymbol{f}_{B_y} =& [\boldsymbol{g}(\boldsymbol{x}^1_{1,},t_{k+1}),\boldsymbol{g}(\boldsymbol{x}^{N_y}_{Q_y,},t_{k+1}),\cdots,\\
&\boldsymbol{g}(\boldsymbol{x}^{(N_x-1)N_y+1}_{1,},t_{k+1}),\boldsymbol{g}(\boldsymbol{x}^{N_xN_y}_{Q_y,},t_{k+1})]^T_{N_xQ_x \times 1},\\
%\end{aligned}
%\end{equation*}
%\begin{equation*}
%    \begin{aligned}
\boldsymbol{f}_{C^0_x} =& [\boldsymbol{0}]_{N_y(N_x-1)Q_y \times 1},\\
\boldsymbol{f}_{C^0_y} =& [\boldsymbol{0}]_{N_x(N_y-1)Q_x \times 1},\\
\boldsymbol{f}_{C^1_x} =& [\boldsymbol{0}]_{N_y(N_x-1)Q_y \times 1},\\
\boldsymbol{f}_{C^1_y} =& [\boldsymbol{0}]_{N_x(N_y-1)Q_x \times 1}.\\
\end{aligned}
\end{equation*}

The coefficient matrices and the right-hand side terms can be rewritten by the rescaling parameters
\begin{equation*}
\begin{aligned} 
\lambda_{i} =&\frac{c}{\mathop{\max}\limits_{1\leq j \leq M_pJ_n}\left|A_{i,j}\right|}, \ A_{i,j} = \lambda_{i}A_{i,j},\  f_{A,i,1} = \lambda_{i}f_{A,i,1},\\
i=&1,\cdots,M_pQ, \\
% &\vdots\\
\lambda_{i} =&\frac{c}{\mathop{\max}\limits_{1\leq j \leq M_pJ_n}\left|B_{x,i,j}\right|},\ B_{x,i,j} = \lambda_{i}B_{x,i,j},\  f_{B_x,i,1} = \lambda_{i}f_{B_x,i,1},\\
i=&1,\cdots,N_yQ_y, \\
\lambda_{i} =&\frac{c}{\mathop{\max}\limits_{1\leq j \leq M_pJ_n}\left|B_{y,i,j}\right|},\ B_{y,i,j} = \lambda_{i}B_{y,i,j},\  f_{B_y,i,1} = \lambda_{i}f_{B_y,i,1},\\
i=&1,\cdots,N_xQ_x, \\
\lambda_{i} =&\frac{c}{\mathop{\max}\limits_{1\leq j \leq M_pJ_n}\left|C^0_{x,i,j}\right|}, \ C^0_{x,i,j} = \lambda_{i}C^0_{x,i,j}, \ f_{C^0_x,i,1} = \lambda_{i}f_{C^0_x,i,1},\\
i=&1,\cdots,N_y(N_x-1)Q_y. \\
\lambda_{i} =&\frac{c}{\mathop{\max}\limits_{1\leq j \leq M_pJ_n}\left|C^0_{y,i,j}\right|}, \ C^0_{y,i,j} = \lambda_{i}C^0_{y,i,j}, \ f_{C^0_y,i,1} = \lambda_{i}f_{C^0_y,i,1},\\
i=&1,\cdots,N_x(N_y-1)Q_x. \\
\lambda_{i} =&\frac{c}{\mathop{\max}\limits_{1\leq j \leq M_pJ_n}\left|C^1_{x,i,j}\right|}, \ C^1_{x,i,j} = \lambda_{i}C^1_{x,i,j}, \ f_{C^1_x,i,1} = \lambda_{i}f_{C^1_x,i,1},\\
i=&1,\cdots,N_y(N_x-1)Q_y. \\
\lambda_{i} =&\frac{c}{\mathop{\max}\limits_{1\leq j \leq M_pJ_n}\left|C^1_{y,i,j}\right|}, \ C^1_{y,i,j} = \lambda_{i}C^1_{y,i,j}, \ f_{C^1_y,i,1} = \lambda_{i}f_{C^1_y,i,1},\\
i=&1,\cdots,N_x(N_y-1)Q_x. \\
\end{aligned}
\end{equation*}

Then, we attain
\begin{equation}
 \mathcal{A}_{t_{k+1}}\boldsymbol{U}_{t_{k+1}}=   \begin{aligned}
\left[\begin{array}{l}
\boldsymbol{A} \\
\boldsymbol{B}_x\\
\boldsymbol{B}_y\\
\boldsymbol{C}^0_x\\
\boldsymbol{C}^0_y\\
\boldsymbol{C}^1_x\\
\boldsymbol{C}^1_y
\end{array}\right]\boldsymbol{U}_{t_{k+1}}=\left[\begin{array}{l}
\boldsymbol{f}_A \\
\boldsymbol{f}_{B_x}\\
\boldsymbol{f}_{B_y}\\
\boldsymbol{f}_{C^0_x}\\
\boldsymbol{f}_{C^0_y}\\
\boldsymbol{f}_{C^1_x}\\
\boldsymbol{f}_{C^1_y}
\end{array}\right]
\end{aligned}=\boldsymbol{f}_{t_{k+1}},
\label{eq:matrix_f}
\end{equation}
related to the loss function \eqref{eq:multiphaseloss}, where $\boldsymbol{U}_{t_{k+1}} = \left[u_{11, t_{k+1}},\cdots,u_{1J_n, t_{k+1}},\cdots,u_{M_p1, t_{k+1}},\cdots,u_{M_pJ_n, t_{k+1}}\right]^T$.
Finally, the linear least-squares method is employed to solve the overdetermined system of equations \eqref{eq:matrix_f} and the optimal $\boldsymbol{U}_{t_{k+1}}=\left(u_{n j, t_{k+1}}\right)^{T}$ is obtained.
\bibliography{journal}% Produces the bibliography via BibTeX.

\end{document}